\documentclass[11pt]{article}

\usepackage[a4paper,margin=2.55cm]{geometry}
\usepackage[T1]{fontenc}
\usepackage[utf8]{inputenc}
\usepackage{lmodern}
\usepackage{microtype}
\usepackage{amsmath,amssymb,amsthm,mathtools}
\usepackage{graphicx}
\usepackage{booktabs}
\usepackage{tabularx}
\usepackage{longtable}
\usepackage{pdflscape}
\usepackage{array}
\usepackage{multirow}
\usepackage{enumitem}
\usepackage{xcolor}
\usepackage{float}
\usepackage[sort&compress,numbers]{natbib}
\usepackage[hidelinks]{hyperref}
\usepackage[nameinlink,noabbrev]{cleveref}
\usepackage{tikz}
\usepackage{standalone}
\usepackage{subcaption}
\usepackage{placeins}

\definecolor{lipblue}{RGB}{35,86,136}
\definecolor{liporange}{RGB}{220,121,43}
\definecolor{lipgray}{RGB}{88,96,105}

\graphicspath{{figures/}}
\allowdisplaybreaks

\newtheorem{definition}{Definition}[section]
\newtheorem{proposition}[definition]{Proposition}
\newtheorem{theorem}[definition]{Theorem}
\newtheorem{corollary}[definition]{Corollary}

\newtheorem{remark}[definition]{Remark}

\floatstyle{ruled}
\newfloat{algorithm}{tbp}{loa}
\floatname{algorithm}{Algorithm}

\usetikzlibrary{arrows.meta,positioning,calc,fit,backgrounds}

\newcolumntype{Y}{>{\raggedright\arraybackslash}X}
\newcolumntype{C}[1]{>{\centering\arraybackslash}p{#1}}
\newcolumntype{R}[1]{>{\raggedleft\arraybackslash}p{#1}}

\title{Robust Pair--Restricted Lipschitz Regression with Global--Local Transfer}

\author{
Sergi Sanjuan, Enrique A. S\'anchez P\'erez, Roger Arnau\\
\small Institut Universitari de Matem\`atica Pura i Aplicada,\\[-1mm]
\small Universitat Polit\`ecnica de Val\`encia, Valencia, Spain\\
\small \texttt{ssansil@upvnet.upv.es, easancpe@mat.upv.es, roger.arnau@upv.es}
}

\date{}
    
\begin{document}
\maketitle
  
\begin{abstract}
We develop a robust pair-restricted Lipschitz regression framework for grouped longitudinal data represented as finite trajectories in a space equipped with a prescribed or learned pseudometric. For a fixed geometry, latent responses and a Lipschitz bound are estimated jointly through a convex optimisation problem with constraints imposed on selected trajectory pairs. We show that the choice of fidelity loss leads to different structural behaviour: Huber regularisation strictly reduces the pair-restricted Lipschitz seminorm whenever the exact interpolant has a finite positive seminorm, whereas $L^1$ fidelity can preserve exact interpolation below a data-dependent threshold. For grouped observations, global and local Lipschitz structures are combined while keeping group identity outside the trajectory geometry. The global component draws on observations from other groups, while the local component uses each group's history to account for group-specific behaviour. We also characterise when group-specific observations can be fitted exactly under a shared Lipschitz bound and establish a regularity bound for the combined predictor. We describe a computational implementation based on constraint generation and illustrate the framework through an application to longitudinal sow-production data.
\end{abstract}



\section{Introduction}
\label{sec:introduction}

Many prediction problems are intrinsically history dependent: the information relevant for predicting the future state of a system is not always contained in a single observation, but in the sequence of states through which that system has evolved. This situation arises in applications such as clinical monitoring, where repeated patient measurements are used to predict future outcomes \citep{liu2017personalized}, and precision livestock farming, where longitudinal farm records are used to forecast production and reproductive performance \citep{gauthier2022}. Two current states that appear similar may therefore have been reached through substantially different recent histories. A natural mathematical response is to treat each finite history as a single object---a trajectory---and to formulate prediction on a space whose elements are such trajectories. Once this viewpoint is adopted, the choice of geometry becomes fundamental: prediction depends not only on the observed responses, but also on how similarities and differences between trajectories are represented.

Trajectory-based metric structures provide a natural framework for representing recent histories as elements of a predictive space. In graph-based settings, \citet{falciani2022} considered finite walks endowed with quasi-pseudometric structures, while \citet{jonard2026} developed a broader construction of metrics and quasi-pseudometrics on hyperspaces of graph states and on spaces of trajectories. Related work by \citet{arnau2025metric} studies weighted metric spaces of graph walks and the use of Lipschitz extensions on these spaces. Graph-based semi-supervised learning also exploits weighted pairwise similarities to regularise label propagation \citep{bozorgnia2026improved}. The present work takes the trajectory-space viewpoint described above as its geometric starting point. However, whereas those constructions focus on defining meaningful geometries from the underlying graph or hyperspace structure, prediction problems also allow the geometry itself to be adapted to the task. Metric-learning methods follow precisely this principle by learning distances from data rather than prescribing them a priori; in regression, for example, \citet{weinberger2007metric} learned a distance directly from predictive performance. Here, trajectories remain the basic predictive objects, while the geometry used to compare them may be selected according to its usefulness for the prediction problem.

Once a trajectory geometry has been specified, Lipschitz regularity provides a way to relate geometric proximity to variation of the fitted response. Classical McShane--Whitney theory guarantees real-valued Lipschitz extensions without increasing the prescribed Lipschitz constant \citep{mcshane1934,whitney1934}, and Lipschitz methods have subsequently been used in classification and regression on general metric spaces \citep{luxburg2004,gottlieb2017}. 
In particular, \citet{gottlieb2017} combine metric-space regression with Lipschitz extension, while alternative extension-based regression constructions include the $p$-average-slope-minimizing approach of \citet{arnau2023measure}. \citet{aziznejad2022} study a regression formulation in which the Lipschitz constant of the fitted function is directly regularised.
A related distinction is that between exact interpolation, which requires the observed responses themselves to satisfy the prescribed Lipschitz geometry, and Lipschitz smoothing, which allows the fitted response to depart from noisy observations. \citet{zaichyk2024} study this smoothing--extension paradigm in a finite-sample regression setting. Building on these approaches, we study how Huber and $L^1$ fidelity affect exact interpolation and Lipschitz regularity in a pair-restricted latent regression problem on trajectory pseudometric spaces, and extend this construction to grouped global--local prediction and sequential updates.

A further difficulty appears when longitudinal observations are organised into identifiable units or groups. This is precisely the case in applications such as clinical monitoring, where repeated measurements belong to individual patients, or livestock forecasting, where observations are naturally associated with individual animals or farms. In such settings, fitting a separate model to each unit limits the possibility of borrowing information from related units, whereas a single pooled model may fail to represent persistent unit-specific behaviour, a tension that is well known in forecasting collections of related time series \citep{monteromanso2021}.
Standard approaches for accounting for group identity include categorical encodings, learned embeddings, and random-effects formulations, particularly when the number of identifiers is large \citep{guo2016,simchoni2021}. In a metric-based Lipschitz framework, however, incorporating group identity directly into the input representation has an additional consequence: the encoding of the groups becomes part of the geometry itself, introducing extra coordinates or requiring a meaningful distance between group identities to be specified or learned together with the trajectory geometry. This can unnecessarily entangle similarity between trajectories with similarity between groups. We therefore propose a group-aware Lipschitz construction that keeps group identity separate from the trajectory geometry while still allowing controlled information sharing across groups.

Sequential prediction introduces an additional requirement: newly observed responses should be incorporated into future predictions without altering decisions already made. Related online and sequential Lipschitz approaches have been studied in \citet{kuzborskij2017,calliess2020}. In the framework developed here, the trajectory geometry and Lipschitz bounds are selected first and then kept fixed, while new observations are incorporated through constrained updates that preserve the previously established regularity structure.

The main mathematical contribution of this work is the construction and structural analysis of a robust, grouped, and sequential Lipschitz regression framework on trajectory pseudometric spaces. For a fixed geometry, the response is estimated through a finite-dimensional regularised problem of the form
$$
\widehat z
\in
\operatorname*{arg\,min}_{z\in\mathbb R^n}
\left\{
\frac{1}{n}\sum_{i=1}^{n}\rho(y_i-z_i)
+
\eta\,Lip_{\mathcal P}(z)
\right\},
$$
where $\mathcal P$ denotes a prescribed family of constrained pairs, introduced formally in \Cref{sec:robust-fit}.
We show that the behaviour of this problem depends strongly on the fidelity loss: Huber regularisation strictly reduces the Lipschitz norm of any nonconstant exact interpolant (whenever its pair-restricted Lipschitz seminorm is finite and positive), whereas $L^1$ fidelity can preserve exact interpolation below a data-dependent threshold.
For grouped observations, a shared Lipschitz bound provides a controlled mechanism for coupling the group-specific fits, and prediction combines information across and within groups through
$$
\widehat Y
=
(1-\alpha)\widehat Y^G
+
\alpha\widehat Y^L,
\qquad
\alpha\in[0,1].
$$
We further establish regularity properties of this global--local predictor and show that the resulting structure can be continued sequentially through fixed-bound updates that remain feasible as new observations are incorporated.

To illustrate how the abstract framework can be instantiated in practice, we include a brief application to longitudinal sow-production data. In this setting, the commercial farms act as the identifiable groups, while recent monthly climatic and production histories are represented as finite trajectories from which the subsequent farrowing rate is predicted. The geometry between trajectories is learned from historical data, so that the abstract notions of trajectory space, group structure, and Lipschitz prediction acquire a concrete interpretation in a longitudinal forecasting problem.

%The remainder of the article is organised as follows. \Cref{sec:background} introduces trajectory spaces, pseudometric geometries, and the Lipschitz-extension results required throughout the paper. \Cref{sec:robust-fit} develops the pair-restricted robust Lipschitz regression problem and analyses its structural behaviour under Huber and $L^1$ fidelity. \Cref{sec:global-local-operators} introduces the grouped global--local construction and the role of the shared local Lipschitz bound. \Cref{sec:update} develops the fixed-bound online continuation and its feasibility properties. \Cref{sec:computational-implementation} describes the computational solution based on constraint generation. \Cref{sec:application} provides the sow-production illustration, and \Cref{sec:conclusions} discusses the scope, limitations, and mathematical implications of the framework.
\section{Theoretical background}
\label{sec:background}

The regression framework developed in the following sections is built on two geometric ingredients. First, recent histories are regarded as points of a space of finite trajectories endowed with a task-adapted pseudometric. Second, once this geometry has been fixed, Lipschitz functions and their extremal extensions provide the mechanism for controlling variation on the observed trajectories and extending fitted responses to new ones. We introduce the trajectory-space construction and the learned pseudometric first, and then recall the Lipschitz extension results required later.

\subsection{Trajectory spaces and pseudometric geometry}
\label{sec:trajectory-spaces}

Let $\mathcal S\subseteq\mathbb R^q$ be a state space, where each state $s=(s^{(1)},\ldots,s^{(q)})$ collects the $q$ observed features used to describe the system at a single time point, and let $k\in\mathbb N$ denote a fixed trajectory length. A finite $k$-step trajectory is an ordered tuple $\tau=(s_1,\ldots,s_k)\in\mathcal S^k$, and we write $\mathcal T_k(\mathcal S):=\mathcal S^k$ for the corresponding trajectory space. In a given prediction problem, only a subset $\mathcal X\subseteq\mathcal T_k(\mathcal S)$ may be admissible. Since the trajectories considered here have finite and fixed length, their state coordinates can be concatenated into a vector $x=(s_1^\top,\ldots,s_k^\top)^\top\in\mathbb R^p$, with $p=qk$. The ordering is part of the representation: each coordinate identifies both a state variable and its position within the trajectory.

This construction is a finite-dimensional specialisation of the trajectory-space viewpoint developed by \citet{jonard2026}, whose starting point is a graph $G=(V,E)$ viewed as a metric or quasi-pseudometric space. In that setting, ordinary graph walks are ordered sequences of vertices, while the broader framework replaces individual vertices by subsets of an underlying metric space and regards a trajectory as a finite or infinite sequence $\phi=(\phi_i)$ with $\phi_i\in\mathcal B(\Omega)$. Metric and quasi-pseudometric structures can then be lifted from the underlying space to hyperspaces of subsets and subsequently to spaces of trajectories. Ordinary graph walks are recovered as a particular case by taking singleton states $\phi_i=\{v_i\}$ with $v_i\in V$. The present work retains this idea of treating an ordered history as a single point of a trajectory space, but uses a different geometric construction: rather than prescribing a trajectory metric through a Hausdorff-type or measure-based hyperspace distance, the geometry is adapted to the prediction task through a finite-dimensional representation of the trajectories.

To define this geometry, we do not assume that all trajectory coordinates contribute equally when two trajectories are compared. We assume that all trajectory coordinates have been scaled to a common scale and associate a non-negative weight $\omega_j$ with each coordinate. For $p$ trajectory coordinates, we define the set of admissible weights as
$$\Delta_p=\left\{\omega\in[0,1]^p:\sum_{j=1}^{p}\omega_j=1\right\}.$$
Given $\omega\in\Delta_p$, the distance between two trajectories $x, x' \in \mathbb R^p$ is defined through the weighted Euclidean (semi-)norm
$$
d_\omega(x,x')
=
\left(
\sum_{j=1}^{p}
\omega_j(x_j-x'_j)^2
\right)^{1/2}.
$$
Thus, each state-feature/position coordinate contributes separately to the comparison of two trajectories, and the weights determine their relative influence on the resulting geometry. 
The condition $\sum_{j=1}^{p}\omega_j=1$ removes the arbitrary freedom to rescale all weights by a common positive factor.

If all weights are strictly positive, $d_\omega$ is a metric on $\mathbb R^p$. Allowing zero weights makes the natural general formulation a pseudometric: distinct trajectories may satisfy $d_\omega(x,x')=0$, although on a particular admissible set $\mathcal X$ the remaining coordinates may still separate all points. The equivalence relation $x\sim_\omega x'$ whenever $d_\omega(x,x')=0$ therefore induces a genuine metric on the quotient space $\mathcal X/{\sim_\omega}$. Retaining the pseudometric formulation is useful because it allows coordinates to receive zero weight without requiring a separate mathematical treatment; the Lipschitz framework introduced below is correspondingly formulated so that zero-distance pairs are handled explicitly.

The geometry is learned through predictive performance rather than from a criterion defined solely on the distances themselves. This type of performance-based weighting is related to metric learning for regression, where distances are selected according to their effect on predictive accuracy \citep{weinberger2007metric}.
Each $\omega\in\Delta_p$ determines a candidate pseudometric $d_\omega$, and if $E_{\mathrm{sel}}(\omega)$ denotes the corresponding prediction error under a prescribed model-selection procedure, the selected geometry satisfies
$$
\widehat\omega
\in
\operatorname*{arg\,min}_{\omega\in\Delta_p}
E_{\mathrm{sel}}(\omega).
$$
Other structural parameters may be selected jointly with $\omega$ when required.
No convexity or uniqueness is assumed for this outer geometry-learning problem; the mathematical results developed in the following sections are conditional on a fixed selected pseudometric $d_{\widehat\omega}$.

\subsection{Lipschitz functions and McShane--Whitney extensions}
\label{sec:lipschitz-framework}

Let $(\mathcal X,d)$ be a (pseudo-)metric space and let $S\subseteq\mathcal X$. A real-valued function $f:S\to\mathbb R$ is said to be Lipschitz if there exists a constant $K\geq 0$ such that
$$
|f(x)-f(x')|\leq Kd(x,x') \qquad\text{ for all } x,x'\in S.
$$
The infimum of all constants $K$ satisfying this inequality is denoted by $Lip(f)$ and is referred to as the Lipschitz (semi-)norm of $f$.

For the trajectory spaces introduced in \Cref{sec:trajectory-spaces}, $d=d_\omega$ is the learned trajectory distance.
Lipschitz regularity then relates variation of the response to the learned geometry: trajectories that are close with respect to $d_\omega$ cannot be assigned arbitrarily different values. The quantity $Lip(f)$ controls the variation of $f$ relative to the trajectory distance. A large value permits substantial differences in $f(x)$ and $f(x')$ even when $d_\omega(x,x')$ is small, whereas a smaller value imposes stronger regularity across nearby trajectories.

The same formulation applies when $d_\omega$ is only a pseudometric. Indeed, if $d_\omega(x,x')=0$, the Lipschitz condition necessarily gives $f(x)=f(x')$. Hence a Lipschitz function is constant on the equivalence classes induced by zero $d_\omega$-distance and can equivalently be regarded as a function on the associated metric quotient.

Suppose now that $A=\{a_1,\ldots,a_n\}\subseteq\mathcal X$ is a finite nonempty set and that $f:A\to\mathbb R$ is $L$-Lipschitz for some $L\geq 0$. The classical McShane--Whitney extension theorem \citep{mcshane1934,whitney1934} gives extremal $L$-Lipschitz extensions of $f$ to $\mathcal X$. With the convention used throughout this article, they are
$$
f^M_L(x)=\max_{a\in A}\{f(a)-L\,d(x,a)\},
\qquad
f^W_L(x)=\min_{a\in A}\{f(a)+L\,d(x,a)\},
\qquad x\in\mathcal X.
$$
Both functions agree with $f$ on $A$ and are $L$-Lipschitz on the whole space $\mathcal X$.
More precisely, \citet[Theorem~4.1.1]{cobzas2019lipschitz} show that if $F:\mathcal X\to\mathbb R$ is any other $L$-Lipschitz extension of the same data, then
$$
f^M_L(x)\leq F(x)\leq f^W_L(x)
\qquad\text{for every }x\in\mathcal X.
$$

Therefore, the difference $f^W_L(x) - f^M_L(x)$ provides a deterministic measure of how strongly the geometry and the available anchors constrain the extension at a query point. Writing $d(x,A)=\min_{a\in A}d(x,a)$, one has
\begin{equation}
\label{eq:rel_width_L}
    0\leq f^W_L(x)-f^M_L(x)\leq 2L\,d(x,A).
\end{equation}

Moreover, the convex combination
$$
E_{\lambda,L}(x)
=
\lambda f^M_L(x)
+
(1-\lambda) f^W_L(x)
, \qquad \lambda\in[0,1],
$$
is also an $L$-Lipschitz extension of $f$. This provides a simple family of predictors between the two extremal extensions without changing the prescribed regularity level.

These results provide the extension mechanism used throughout the paper once a trajectory geometry and compatible anchor values have been fixed. The next section addresses the complementary fitting problem: how to estimate those anchor values and their Lipschitz regularity from noisy observations.
\section{Robust Lipschitz regression}
\label{sec:robust-fit}

The extension framework of \Cref{sec:lipschitz-framework} assumes that compatible anchor values and a Lipschitz bound are already available. In real prediction problems, however, observed responses may contain measurement noise, atypical observations, or short-term fluctuations that are not well represented by the underlying trajectory geometry. Enforcing exact interpolation in such cases can force nearby trajectories to reproduce large response differences and may therefore require an unnecessarily large Lipschitz bound. Rather than treating every observed response as an exact geometric constraint, we replace the observations by latent fitted values that remain close to the data while being regularised by their Lipschitz variation. This leads to a finite-dimensional variational problem that provides a robust interface between the learned trajectory geometry and the extension machinery developed above. We first establish its basic convexity and existence properties and then analyse how its behaviour depends on the choice between Huber and $L^1$ fidelity.

\subsection{Pair-restricted variational formulation}
\label{sec:pair-restricted-fit}

Fix the learned trajectory geometry $d_\omega$ and consider a finite repository $\mathcal R=\{(x_i,y_i,g_i)\}_{i=1}^{n}$, where $x_i\in\mathcal X$ is a trajectory, $y_i\in\mathbb R$ is its observed response, and $g_i\in\mathcal G$ identifies the group from which it was obtained. The group labels are retained for the grouped construction developed in \Cref{sec:global-local-operators}; in the present section they enter only through the choice of the constrained pair family $\mathcal P$. Rather than extending the observed values $y_i$ directly, we associate each trajectory $x_i$ with a latent value $z_i$. The latent values are fitted to remain close to the observations while controlling their variation with respect to the learned trajectory geometry.

Let $\mathcal P\subseteq\{(i,j):1\leq i<j\leq n\}$ denote the family of trajectory pairs on which Lipschitz compatibility is imposed, and let $\rho:\mathbb{R}\to[0,\infty)$ be a convex loss function.
We estimate the latent values $z=(z_1,\ldots,z_n)$ together with a Lipschitz bound $\ell\geq0$ by solving
\begin{equation}
\label{pr:latent-fit}
\begin{aligned}
(\widehat z,\widehat\ell)
\in
\operatorname*{arg\,min}_{z\in\mathbb R^n,\,\ell\geq0}
\quad&
\frac{1}{n}\sum_{i=1}^{n}\rho(y_i-z_i)
+\eta\ell,
\\
\text{subject to}\quad&
|z_i-z_j|
\leq
\ell\,d_\omega(x_i,x_j),
\qquad
(i,j)\in\mathcal P .
\end{aligned}
\end{equation}
Here the coefficient $\eta>0$ controls the trade-off between fidelity to the observations and Lipschitz regularity. Larger values penalise the Lipschitz bound more strongly and favour a more regular latent response, potentially at the cost of a larger discrepancy from the observations. Smaller values give greater weight to data fidelity and allow more variation with respect to the trajectory geometry.

The Lipschitz regularity of a latent response vector can be expressed directly. For $z\in\mathbb R^n$, define its extended pair-restricted Lipschitz (semi-)norm by
$$
Lip_{\mathcal P}(z)=
\begin{cases}
\displaystyle \max_{\substack{(i,j)\in\mathcal P\\ d_\omega(x_i,x_j)>0}}
\frac{|z_i-z_j|}{d_\omega(x_i,x_j)},
& \text{if } z_i=z_j \text{ whenever } (i,j)\in\mathcal P \text{ and } d_\omega(x_i,x_j)=0,\\[3mm]
+\infty,
& \text{otherwise},
\end{cases}
$$
with the convention that the maximum over an empty set is zero.
For a fixed $z$, this is precisely the smallest admissible value of $\ell$ in Problem~\ref{pr:latent-fit}. This leads to the equivalent variational formulation
\begin{equation}
\label{eq:variational-fit}
\widehat z
\in
\operatorname*{arg\,min}_{z\in\mathbb R^n}
\left\{
\frac{1}{n}\sum_{i=1}^{n}\rho(y_i-z_i)
+
\eta\,Lip_{\mathcal P}(z)
\right\}.
\end{equation}

The variational structure in \Cref{eq:variational-fit} has several precedents in Lipschitz-regularised learning. In metric-space classification, \citet{luxburg2004} penalise the empirical loss of the fitted values by their finite-sample Lipschitz constant. In regression, \citet{aziznejad2022} study a continuous-domain univariate formulation in which a data-fidelity term is penalised directly by the Lipschitz constant of the fitted function, while Lipschitz constraints have also been considered in other regression settings as in Lipschitz-constrained isotonic regression studied by \citet{yeganova2009}. On graphs, \citet{kyng2015} study Lipschitz learning and smooth extensions, including $L^1$ regularisation of observed values and edge weights. \citet{zaichyk2024} formulate Lipschitz smoothing with compatibility constraints indexed by the edges of a prescribed graph, with the complete graph recovering the full pairwise constraint family.
When $\mathcal P$ contains all trajectory pairs, $Lip_{\mathcal P}(z)$ recovers the corresponding finite-sample Lipschitz regulariser. The formulation considered here differs in operating on a learned trajectory (pseudo)metric, allowing a prescribed family of constrained pairs and using the resulting latent values as anchors for the subsequent global--local extension.

Two convex fidelity losses are particularly useful for the structural analysis below. The first is the absolute loss (or $L^1$ loss), $\rho_{L^1}(r)=|r|$.
The second is the Huber loss \cite{huber1964},
\begin{equation*}
\rho_\delta(r)=
\begin{cases}
\frac{1}{2}r^2,
& |r|\leq \delta,\\[1mm]
\delta |r|-\frac{1}{2}\delta^2,
& |r|>\delta,
\end{cases}
\qquad \delta>0.
\end{equation*}
The $L^1$ loss grows linearly with the residual magnitude and is therefore less sensitive to large discrepancies than the squared loss.
On the other hand, the Huber loss combines the quadratic behaviour of the squared loss for small residuals with linear growth in the tails, thereby reducing the influence of large discrepancies while retaining a smooth quadratic regime around the origin. The parameter $\delta$ determines the transition between these two regimes.

Both the $L^1$ and Huber losses are finite, continuous, convex, and coercive. Recall that a function $\rho:\mathbb{R}\to[0,\infty)$ is coercive if $\rho(r)\to+\infty$ as $|r|\to\infty$. They therefore fall within the class of loss functions covered by the following result, which establishes the basic well-posedness properties of the latent Lipschitz fitting problem.

\begin{proposition}%[Convexity and existence of the latent Lipschitz fit]
\label{prop:latent-fit}
Let $\rho:\mathbb{R}\to[0,\infty)$ be a convex and coercive loss function, and let $\eta>0$. Fix a finite nonempty repository $\mathcal R$, a trajectory geometry $d_\omega$, and a family of constrained pairs $\mathcal P$. Then Problem~\ref{pr:latent-fit} is a finite-dimensional convex optimisation problem with a nonempty feasible set and at least one minimiser. Moreover, every minimiser satisfies
$$
\widehat\ell=Lip_{\mathcal P}(\widehat z),
$$
so Problem~\ref{pr:latent-fit} is equivalently given by the variational formulation in \cref{eq:variational-fit}. If, in addition, $\rho$ is strictly convex, then the minimiser is unique.
\end{proposition}
\begin{proof}
Since the trajectory geometry $d_\omega$ is fixed, for every $(i,j)\in\mathcal P$, the Lipschitz constraint $|z_i-z_j|\leq\ell\,d_\omega(x_i,x_j)$ is equivalent to the two affine inequalities
$$
z_i-z_j-\ell\,d_\omega(x_i,x_j)\leq0,
\qquad
z_j-z_i-\ell\,d_\omega(x_i,x_j)\leq0.
$$
Hence the feasible set is closed and convex. It is nonempty because $z_i=c$ for all $i$ and $\ell=0$ is feasible for every $c\in\mathbb R$. The objective
$$
F(z,\ell)=\frac{1}{n}\sum_{i=1}^{n}\rho(y_i-z_i)+\eta\ell
$$
is continuous and convex because the loss function $\rho$ is continuous and convex and $\eta \ell$ is linear.

It remains to show that the infimum is attained. Choose any feasible point and let $C$ be its finite objective value. On the feasible sublevel set with objective at most $C$, $0\leq\ell\leq C/\eta$. Moreover,
$$
\rho(y_i-z_i)\leq nC
$$
for every $i$. Since $\rho(r)\to\infty$ as $|r|\to\infty$ by coercivity, each $z_i$ is bounded on this sublevel set. The sublevel set is therefore closed and bounded in finite dimension, hence compact. By the Weierstrass theorem, $F(z,\ell)$ attains its minimum and Problem~\ref{pr:latent-fit} admits at least one minimiser.

For any fixed $z$ with $Lip_{\mathcal P}(z)<\infty$, the smallest feasible value of $\ell$ is precisely $Lip_{\mathcal P}(z)$.
Since $\eta>0$, a minimiser cannot have $\ell>Lip_{\mathcal P}(z)$, as decreasing $\ell$ to $Lip_{\mathcal P}(z)$ would preserve feasibility and strictly decrease the objective. Hence every minimiser satisfies $\widehat\ell=Lip_{\mathcal P}(\widehat z)$, which yields the equivalent variational formulation in \cref{eq:variational-fit}.

Finally, if $\rho$ is strictly convex, then the fidelity term is strictly convex in $z$. Since $Lip_{\mathcal P}$ is convex, the variational objective is strictly convex on its effective domain and therefore has at most one minimiser $\widehat z$. The identity $\widehat\ell=Lip_{\mathcal P}(\widehat z)$ then also determines $\widehat\ell$ uniquely.
\end{proof}

\begin{remark}\label{rem:l1-nonuniqueness}
Without strict convexity of the loss, uniqueness is not guaranteed. For example, consider the absolute loss $\rho(r)=|r|$, two observations $y_1=0$ and $y_2=1$, a single constrained pair with $d_\omega(x_1,x_2)=1$, and $\eta \geq 1/2$. Since feasibility requires $\ell\geq|z_1-z_2|$,
$$
F(z,\ell)
=
\frac12\bigl(|z_1|+|1-z_2|\bigr)+\eta\ell
\geq
\frac12\bigl(|z_1|+|z_1-z_2|+|1-z_2|\bigr)
\geq
\frac12.
$$
Every choice $0 \leq z_1 = z_2 \leq 1$ and $\ell=0$ attains this value, so the minimiser is not unique. In particular, neither the $L^1$ loss nor the Huber loss is strictly convex on the whole real line, and uniqueness therefore cannot be asserted in general for the fits considered here.
\end{remark}

The preceding result establishes well-posedness but does not determine whether the fitted response coincides with the observations. This question depends fundamentally on the fidelity loss and is analysed next.

\subsection{Structural properties under Huber and \texorpdfstring{$L^1$}{L1} fidelity}
\label{sec:fidelity-structure}

The Huber formulation has an additional structural consequence that clarifies the role of the Lipschitz penalty. Whenever the observed response has a finite and positive pair-restricted Lipschitz norm, any positive Huber--Lipschitz regularisation forces the fitted latent response away from exact interpolation and yields a strictly smaller pair-restricted Lipschitz norm. The following result makes this property precise.

\begin{theorem}%[Strict Lipschitz shrinkage under Huber–Lipschitz regularisation]
\label{thm:huber-strict-shrinkage}
Assume that the Huber loss $\rho_\delta$, with $\delta>0$, is used in \cref{eq:variational-fit} for some $\eta>0$, and that the observed response vector $y=(y_1,\ldots,y_n)$ satisfies $0<Lip_{\mathcal P}(y)<\infty$.
Then the exact interpolant $z=y$ is not a minimiser of the latent problem. Moreover, every minimiser $z^\ast$ satisfies
$$
Lip_{\mathcal P}(z^\ast)<Lip_{\mathcal P}(y).
$$
\end{theorem}
\begin{proof}
Fix any $c\in\mathbb R$ and consider
$$
z^{(t)}=(1-t)y+t\,c\mathbf 1,\qquad 0<t<1,
$$
where $\mathbf 1=(1,\ldots,1)$. By construction,
$$
Lip_{\mathcal P}\bigl(z^{(t)}\bigr)=(1-t)Lip_{\mathcal P}(y).
$$
Moreover, $y_i-z_i^{(t)}=t(y_i-c)$. Because the repository is finite, for sufficiently small $t>0$ all these residuals satisfy $|t(y_i-c)|\leq\delta$ and therefore lie in the quadratic region of the Huber loss. Hence
$$
\frac{1}{n}\sum_{i=1}^{n}\rho_\delta\bigl(y_i-z_i^{(t)}\bigr)
=
\frac{t^2}{2n}\sum_{i=1}^{n}(y_i-c)^2.
$$
Let $J$ be the regularised objective in \cref{eq:variational-fit}. Then
$$
J\bigl(z^{(t)}\bigr)-J(y)
=
\frac{t^2}{2n}\sum_{i=1}^{n}(y_i-c)^2-\eta t Lip_{\mathcal P}(y)
=
t \left( \frac{t}{2n}\sum_{i=1}^{n}(y_i-c)^2-\eta Lip_{\mathcal P}(y) \right).
$$
The term in parentheses converges to $-\eta Lip_{\mathcal P}(y)<0$ as $t\to0^+$.
Therefore, for all sufficiently small $t>0$, it holds that $J\bigl(z^{(t)}\bigr)<J(y)$, so that $z=y$ cannot be a minimiser.

Now let $z^\ast$ be any minimiser. Since $\rho_\delta\geq0$,
$$
\eta\,Lip_{\mathcal P}(z^\ast)
\leq
J(z^\ast)
\leq
J\bigl(z^{(t)}\bigr)
<
J(y)
=
\eta\,Lip_{\mathcal P}(y)
$$
for sufficiently small $t>0$. Since $\eta>0$, it follows that $Lip_{\mathcal P}(z^\ast)<Lip_{\mathcal P}(y)$.
\end{proof}

The argument is not specific to the Huber loss: the same strict-shrinkage conclusion holds for any nonnegative fidelity loss satisfying $\rho(0)=0$ and $\rho(r)=\mathcal{O}(|r|^p)$ as $r\to0$ with $p>1$, provided that the corresponding minimisation problem admits a minimiser. The condition $Lip_{\mathcal P}(y)<\infty$ requires the observed responses to be compatible with every constrained pair at zero trajectory distance; otherwise exact interpolation is not feasible.

Unlike the Huber case in \cref{thm:huber-strict-shrinkage}, exact interpolation under the absolute loss may persist for positive regularisation strengths. Because the $L^1$ fidelity term increases linearly with the displacement $\|y-z\|_1$, a reduction in the Lipschitz norm is beneficial only when the corresponding decrease in the regularisation penalty compensates for this fidelity cost. The persistence of an unregularised solution under positive regularisation is related to the broader theory of exact regularisation for convex programs \citep{friedlander2008exact}. Related threshold phenomena also arise in $TV$--$L^1$ regularisation; see, for example, \citet[Lemma~5.5]{chan2005aspects} and \citet[Theorem~3.2 and Corollary~4.2]{choksi2011anisotropic}. The following result establishes a finite-dimensional exact-interpolation characterisation specific to the pair-restricted Lipschitz seminorm used here.

\begin{proposition}%[Exact interpolation under $L^1$ fidelity]
\label{prop:l1-exact-interpolation}
Assume that the absolute loss is used in \cref{eq:variational-fit}, with $\eta>0$, and that $Lip_{\mathcal P}(y)<\infty$.
Define
$$
\kappa(y)
=
\sup_{\substack{z\neq y\\ Lip_{\mathcal P}(z)<Lip_{\mathcal P}(y)}}
\frac{
Lip_{\mathcal P}(y)-Lip_{\mathcal P}(z)
}{
\|y-z\|_1
},
$$
with the convention that the supremum of an empty family is zero. Then
$$
y\in
\operatorname*{arg\,min}_{z\in\mathbb R^n}
\left\{
\frac1n\|y-z\|_1
+
\eta\,Lip_{\mathcal P}(z)
\right\}
\quad\Longleftrightarrow\quad
\eta\,\kappa(y)\leq\frac1n.
$$
If
$
\eta\,\kappa(y)<\frac{1}{n},
$
then $y$ is the unique minimiser.
\end{proposition}
\begin{proof}
Let $J(z)=\frac1n\|y-z\|_1+\eta\,Lip_{\mathcal P}(z)$.
Since $J(y)=\eta\,Lip_{\mathcal P}(y)$, for every $z$ with $Lip_{\mathcal P}(z)<\infty$ we have
$$
J(z)-J(y)
=
\frac1n\|y-z\|_1
-
\eta\left(
Lip_{\mathcal P}(y)-Lip_{\mathcal P}(z)
\right).
$$

If $Lip_{\mathcal P}(z)\geq Lip_{\mathcal P}(y)$ and $z\neq y$, then
$$
J(z)-J(y)
\geq
\frac1n\|y-z\|_1
>
0.
$$
Hence only vectors satisfying $Lip_{\mathcal P}(z)<Lip_{\mathcal P}(y)$ can improve upon the exact interpolant.
In this case,
$$
J(z)-J(y)
=
\|y-z\|_1
\left(
\frac{1}{n}
-
\eta
\frac{
Lip_{\mathcal P}(y)-Lip_{\mathcal P}(z)
}{
\|y-z\|_1
}
\right).
$$
Therefore, $J(z)\geq J(y)$ for every $z$ if and only if
$$
\eta
\frac{
Lip_{\mathcal P}(y)-Lip_{\mathcal P}(z)
}{
\|y-z\|_1
}
\leq
\frac{1}{n}
$$
for every $z\neq y$ with $Lip_{\mathcal P}(z)<Lip_{\mathcal P}(y)$. By the definition of $\kappa(y)$, this is equivalent to
$
\eta\,\kappa(y)\leq\frac{1}{n}$.

If the inequality is strict, then every $z\neq y$ with smaller Lipschitz norm satisfies $J(z)>J(y)$, while vectors with Lipschitz norm at least $Lip_{\mathcal P}(y)$ already satisfy the same strict inequality. Hence $y$ is the unique minimiser.
\end{proof}

The exact threshold in \cref{prop:l1-exact-interpolation} is determined by $\kappa(y)$, which measures the largest possible decrease in the pair-restricted Lipschitz norm per unit of $L^1$ displacement from the observed response. Although this quantity gives an exact characterisation, it need not be easy to evaluate directly. A simpler sufficient condition can be obtained from the smallest positive distance among the constrained trajectory pairs. If at least one pair $(i,j)\in\mathcal P$ satisfies $d_\omega(x_i,x_j)>0$, define
$$
d_{\min}
=
\min_{\substack{(i,j)\in\mathcal P\\ d_\omega(x_i,x_j)>0}}
d_\omega(x_i,x_j).
$$
Then $\kappa(y)\leq 1/d_{\min}$, which yields the following directly checkable interpolation regime.

\begin{corollary}
\label{cor:l1-distance-threshold}
Assume that $\mathcal P$ contains at least one pair $(i,j)$ with $d_\omega(x_i,x_j)>0$ and suppose that $Lip_{\mathcal P}(y)<\infty$. If $\eta\leq\frac{d_{\min}}{n}$, then $y$ is a minimiser of the $L^1$-regularised problem. If $\eta<\frac{d_{\min}}{n}$, then $y$ is the unique minimiser.
\end{corollary}
\begin{proof}
For every $z$ with $Lip_{\mathcal P}(z)<\infty$ we have that
$
Lip_{\mathcal P}(y)
\leq
Lip_{\mathcal P}(z)
+
Lip_{\mathcal P}(y-z).
$
Moreover, for every constrained pair $(i,j)$ with positive distance,
$$
|(y_i-z_i)-(y_j-z_j)|
\leq
|y_i-z_i|+|y_j-z_j|
\leq
\|y-z\|_1,
$$
and therefore
$
Lip_{\mathcal P}(y-z)
\leq
\frac{\|y-z\|_1}{d_{\min}}.
$
Hence
$$
\frac{Lip_{\mathcal P}(y)-Lip_{\mathcal P}(z)}
{\|y-z\|_1}
\leq
\frac1{d_{\min}}
$$
whenever $z\neq y$ and $Lip_{\mathcal P}(z)<Lip_{\mathcal P}(y)$. Taking the supremum over such $z$ yields
$
\kappa(y)\leq\frac1{d_{\min}}.
$

Therefore, if $\eta\leq d_{\min}/n$, then
$
\eta\,\kappa(y)\leq\frac1n,
$
and \cref{prop:l1-exact-interpolation} shows that $y$ is a minimiser. If $\eta<d_{\min}/n$, then $\eta\,\kappa(y)<1/n$, so the same proposition gives uniqueness.
\end{proof}

The strict inequality in \cref{cor:l1-distance-threshold} cannot in general be relaxed when uniqueness is required. In the two-point example considered above, with $y_1=0$, $y_2=1$, $d_\omega(x_1,x_2)=1$, and $\eta=1/2=d_{\min}/n$, the exact interpolant $z=y$ is a minimiser, but every constant vector $z_1=z_2=c$, $c\in[0,1]$, attains the same objective value. Thus, equality guarantees interpolation but does not guarantee uniqueness.

Once a minimiser $(\widehat z,\widehat\ell)$ has been obtained, the fitted values $\widehat z_i$ provide the anchors for the Lipschitz extensions developed in the grouped construction of \Cref{sec:global-local-operators}. We refer to the trajectories together with these fitted values as the latent repository.

\begin{figure}[!t]
\centering
\begin{subfigure}[t]{\linewidth}
\centering
\includegraphics[width=0.94\linewidth]{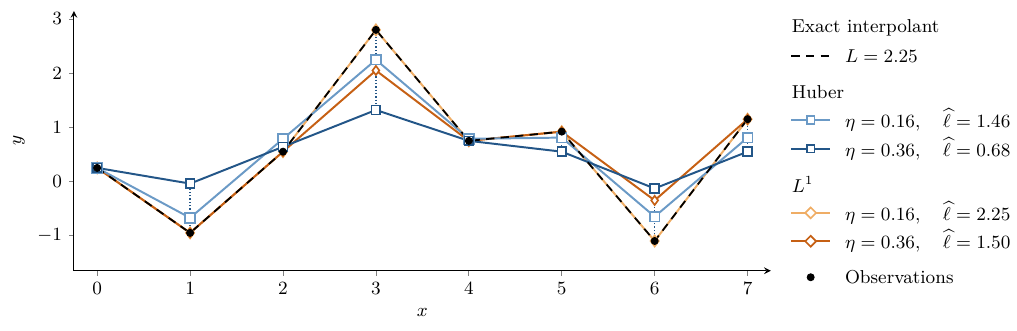}
\caption{Exact interpolation and robust latent fits.}
\label{fig:toy-losses}
\end{subfigure}
\vspace{0.8\baselineskip}
\begin{subfigure}[t]{\linewidth}
\centering
\includegraphics[width=0.94\linewidth]{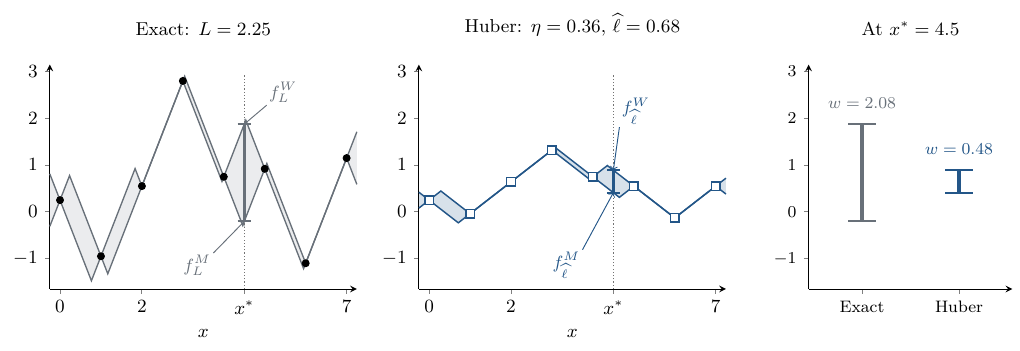}
\caption{McShane--Whitney admissible extension intervals.}
\label{fig:toy-extensions}
\end{subfigure}
\caption{Synthetic illustration of Lipschitz regression on $(\mathbb R,|\cdot|)$ using eight observations and all pairwise constraints.
(a) Exact interpolation and latent fits under Huber and $L^1$ fidelity for $\eta\in\{0.16,0.36\}$.
(b) McShane--Whitney extension envelopes for the exact observations and the Huber fit at $\eta=0.36$. The shaded regions represent admissible Lipschitz extensions satisfying $f_K^M(x)\leq F(x)\leq f_K^W(x)$ for the corresponding Lipschitz bound $K$.}
\label{fig:toy-lipschitz-regression}
\end{figure}

The different behaviour of Huber and $L^1$ fidelity is illustrated in \Cref{fig:toy-lipschitz-regression} through a simple one-dimensional example. For $\eta=0.16$, the $L^1$ fit coincides with the exact interpolant, whereas Huber regularisation already modifies the latent values and strictly reduces the Lipschitz bound, consistently with \Cref{thm:huber-strict-shrinkage,prop:l1-exact-interpolation}. The effect on the McShane--Whitney extensions is also illustrated: at $x^*=4.5$, the admissible interval width decreases from $2.08$ under exact interpolation to $0.48$ for the Huber fit with $\eta=0.36$. This reduction reflects the geometry of the fitted anchors and their Lipschitz bounds, although a smaller Lipschitz bound does not generally guarantee narrower extension intervals at every point.

\section{Global--local regression for grouped observations}
\label{sec:global-local-operators}

As already discussed in the Introduction, many longitudinal prediction problems involve repeated observations from identifiable units, such as patients, farms, machines, or monitoring sites. Observations from different units may share common predictive patterns, while each unit may also exhibit persistent behaviour that is not fully captured by the trajectory coordinates. A purely local model preserves this unit-specific structure but may use only a limited amount of information, whereas a fully pooled model can exploit the complete repository at the cost of overlooking systematic differences between units. The pair-restricted formulation of \Cref{sec:robust-fit} provides a natural way to balance these two sources of information without making group identity part of the trajectory geometry. We use this flexibility to construct two complementary Lipschitz structures: a global fit that couples the complete repository and a local fit that preserves within-group structure under a shared regularity bound. These two components are subsequently extended through the McShane--Whitney operators and combined into a global--local predictor whose Lipschitz regularity remains explicitly controlled.

\subsection{Global and local Lipschitz structures}
\label{sec:global-local-structures}

Repeated observations from identifiable groups create a second modelling problem. Using only the history of the target group preserves group-specific behaviour but limits the amount of information available for prediction, whereas pooling all observations increases the effective repository at the risk of ignoring persistent differences between groups.

A common approach in supervised learning is to introduce group identity as a categorical predictor. One-hot encoding represents each category with a separate binary coordinate, although its dimension grows with the number of identifiers. Another widely used strategy is to represent categories through learned embeddings; for example, neural-network models can learn low-dimensional vector representations of categorical identities jointly with the prediction task \citep{guo2016}.

Such representations could also be incorporated into the trajectory vectors considered here. However, doing so would make group identity part of the learned geometry itself: one-hot encoding would add categorical coordinates to $d_\omega$, whereas embeddings would require learning an additional latent representation of the groups. We instead keep the trajectory metric focused on similarities between observed histories, while group identity determines how the Lipschitz fit is constrained and which historical trajectories can be used when predicting a given group.

We introduce two complementary latent fits. For the global component, Lipschitz regularity is imposed over the complete historical repository,
$$
\mathcal P^G=\{(i,j):1\leq i<j\leq n\},
$$
whereas the local component only constrains trajectories belonging to the same group,
$$
\mathcal P^L=\{(i,j):1\leq i<j\leq n,\ g_i=g_j\}.
$$
Solving Problem~\ref{pr:latent-fit} with $\mathcal P^G$ and $\mathcal P^L$ gives $(\widehat z^G,\widehat\ell_G)$ and $(\widehat z^L,\widehat\ell_L)$, respectively. Both fits use the complete historical repository but differ in their constrained pair families.

For a target group $g\in\mathcal G$, define the index sets used to select the anchors for local and global prediction by
$$
\mathcal I_g^L=\{i:g_i=g\},
\qquad
\mathcal I_g^G=\{i:g_i\neq g\}.
$$
The within-group restricted pair family is
$$
\mathcal P_g^L
=
\left\{(i,j)\in\mathcal P^L:i,j\in\mathcal I_g^L\right\}.
$$
Thus,
$
\mathcal P^L
=
\bigcup_{g\in\mathcal G}\mathcal P_g^L,
$
and no pair in $\mathcal P^L$ connects two different groups. Consequently, once the common local Lipschitz bound is fixed, the optimisation over the local latent values separates across groups. The groups remain coupled only through the shared bound $\widehat\ell_L$, which is estimated jointly in the local fit.
This shared bound is therefore the only coupling mechanism between the groupwise local fits, and its effect on the fitted responses can be characterised explicitly.

\subsection{Shared local regularity and global--local prediction}
\label{sec:shared-local-prediction}

For a given group $g$, let $y_g=(y_i)_{i\in\mathcal I_g^L}$ and $z_g=(z_i)_{i\in\mathcal I_g^L}$ be the observed-response and local latent-response vectors, respectively. Whenever $Lip_{\mathcal P_g^L}(y_g)<\infty$, this quantity is the smallest Lipschitz bound under which the observations from group $g$ can be interpolated exactly. The shared local bound therefore has a direct structural interpretation: groups may require different levels of within-group regularity to interpolate their observations, but the local fit selects a single common bound for all of them. Both the $L^1$ and Huber losses considered in this work satisfy the loss assumptions of the following result.

\begin{theorem}%[Groupwise threshold and saturation under a shared Lipschitz bound]
\label{thm:shared-local-threshold}
Assume that the local latent fit uses a convex loss $\rho:\mathbb R\to[0,\infty)$ satisfying $\rho(0)=0$ and $\rho(r)>0$ for every $r\neq0$, and let $(\widehat z^L,\widehat\ell_L)$ be any minimiser of Problem~\ref{pr:latent-fit} with $\mathcal P=\mathcal P^L$. For each group $g\in\mathcal G$ such that $Lip_{\mathcal P_g^L}(y_g)<\infty$, the following statements hold:
\begin{enumerate}
    \item $\widehat z_g^L=y_g \Longleftrightarrow Lip_{\mathcal P_g^L}(y_g)\leq\widehat\ell_L,$
    \item $Lip_{\mathcal P_g^L}(\widehat z_g^L) = \min\left\{Lip_{\mathcal P_g^L}(y_g),\widehat\ell_L\right\}.$
\end{enumerate}
\end{theorem}
\begin{proof}
Fix $\widehat\ell_L$. Since the local constraints only connect trajectories belonging to the same group, conditional on this common bound the optimisation over the latent values separates across groups. Hence, for each $g\in\mathcal G$, $\widehat z_g^L$ minimises the fidelity loss subject to
$$
Lip_{\mathcal P_g^L}(z_g)\leq\widehat\ell_L.
$$

For (1), suppose first that $Lip_{\mathcal P_g^L}(y_g)\leq\widehat\ell_L$. Then $y_g$ is feasible and has zero fidelity loss. Since $\rho(r)>0$ for every $r\neq0$, no other vector can attain zero fidelity loss, and therefore $\widehat z_g^L=y_g$. Conversely, if $\widehat z_g^L=y_g$, feasibility immediately gives $Lip_{\mathcal P_g^L}(y_g)\leq\widehat\ell_L$.

For (2), consider first the case $Lip_{\mathcal P_g^L}(y_g)\leq\widehat\ell_L$. By (1), $\widehat z_g^L=y_g$, and hence
$$
Lip_{\mathcal P_g^L}(\widehat z_g^L)
=
Lip_{\mathcal P_g^L}(y_g)
=
\min\left\{
Lip_{\mathcal P_g^L}(y_g),
\widehat\ell_L
\right\}.
$$

It remains to consider the case $Lip_{\mathcal P_g^L}(y_g)>\widehat\ell_L$. By (1), $\widehat z_g^L\neq y_g$. Suppose, for contradiction, that $Lip_{\mathcal P_g^L}(\widehat z_g^L)<\widehat\ell_L$. For sufficiently small $t>0$, consider
$$
\widetilde z_g=(1-t)\widehat z_g^L+t y_g.
$$
By convexity of $Lip_{\mathcal P_g^L}$,
$$
Lip_{\mathcal P_g^L}(\widetilde z_g)
\leq
(1-t)Lip_{\mathcal P_g^L}(\widehat z_g^L)
+
t\,Lip_{\mathcal P_g^L}(y_g)
<
\widehat\ell_L,
$$
for all sufficiently small $t>0$, so $\widetilde z_g$ remains feasible. On the other hand,
$$
y_g-\widetilde z_g=(1-t)(y_g-\widehat z_g^L).
$$
Thus every nonzero residual is multiplied by the factor $1-t\in(0,1)$. By convexity of $\rho$ and since $\rho(0)=0$, for every nonzero residual $r$,
$$
\rho((1-t)r)
\leq
(1-t)\rho(r)+t\rho(0)
=
(1-t)\rho(r)
<
\rho(r).
$$
Since $\widehat z_g^L\neq y_g$, at least one residual is nonzero, and therefore the total fidelity loss at $\widetilde z_g$ is strictly smaller than at $\widehat z_g^L$. This contradicts the optimality of $\widehat z_g^L$. Hence, $Lip_{\mathcal P_g^L}(\widehat z_g^L)=\widehat\ell_L$.
Since $Lip_{\mathcal P_g^L}(y_g)>\widehat\ell_L$ in this case, it follows that
$
Lip_{\mathcal P_g^L}(\widehat z_g^L)
=
\min\left\{
Lip_{\mathcal P_g^L}(y_g),
\widehat\ell_L
\right\}.
$
\end{proof}

The theorem shows that the shared local Lipschitz bound acts as an exact threshold across groups. If the observed within-group Lipschitz norm does not exceed $\widehat\ell_L$, the observations of that group are retained without latent adjustment. If the observed norm exceeds $\widehat\ell_L$, the fitted latent response is regularised exactly up to the common bound.
Thus, $\widehat\ell_L$ is not only a common upper bound on within-group variation: it separates groups that can be interpolated exactly from those whose latent responses must be adjusted and determines the regularity level attained by the latter.

The fitting and prediction scopes are different. Although the global latent fit is estimated jointly from the complete historical repository using $\mathcal P^G$, its extension for a target group $g$ uses only the fitted anchors indexed by $\mathcal I_g^G$, whereas the local extension uses the fitted anchors indexed by $\mathcal I_g^L$.

For the prediction formulas below, we restrict attention to target groups $g$ for which both anchor sets $I_g^G$ and $I_g^L$ are nonempty. The restrictions of the fitted global and local repositories to these sets remain $\widehat\ell_G$- and $\widehat\ell_L$-Lipschitz, respectively, so the corresponding McShane--Whitney extensions are well defined. If either anchor set is empty, an additional prediction rule would be required and is outside the present formulation.

Using the global latent fit, the McShane--Whitney extensions for group $g$ at time $t$ are
\begin{equation*}
M^G_{g,t}
=
\max_{i\in\mathcal I_g^G}
\left\{
\widehat z_i^G
-
\widehat\ell_G\,d_\omega(x_{g,t},x_i)
\right\},
\qquad
W^G_{g,t}
=
\min_{i\in\mathcal I_g^G}
\left\{
\widehat z_i^G
+
\widehat\ell_G\,d_\omega(x_{g,t},x_i)
\right\}.
\end{equation*}
The corresponding global prediction is
\begin{equation*}
\widehat Y^G_{g,t}
=
\lambda_G M^G_{g,t}
+
(1-\lambda_G) W^G_{g,t},
\qquad
\lambda_G\in[0,1].
\end{equation*}
The local prediction is constructed analogously from the local latent fit and the same-group index set $\mathcal I_g^L$, yielding $\widehat Y^L_{g,t}$ with mixing parameter $\lambda_L\in[0,1]$.

Finally, the global and local predictions are combined as
\begin{equation*}
\widehat Y_{g,t}
=
\alpha\widehat Y^L_{g,t}
+
(1-\alpha)\widehat Y^G_{g,t},
\qquad
\alpha\in[0,1],
\end{equation*}
where the coefficient $\alpha$ controls the balance between information transferred from the remaining groups and information specific to the target group.

For a fixed target group $g$, repository, trajectory geometry, latent values, Lipschitz bounds, and mixing coefficients, view the global and local predictors as functions of the query trajectory $x$. Since the McShane and Whitney extensions in each branch are Lipschitz with bounds $\widehat\ell_G$ and $\widehat\ell_L$, respectively, their convex combinations satisfy
$$
Lip(\widehat Y_g^G)\leq\widehat\ell_G,
\qquad
Lip(\widehat Y_g^L)\leq\widehat\ell_L.
$$
Consequently, the global--local predictor satisfies
$$
Lip(\widehat Y_g)
\leq
\alpha\widehat\ell_L + (1-\alpha)\widehat\ell_G.
$$

The global--local construction therefore combines cross-group transfer and group-specific information while retaining an explicit Lipschitz regularity bound. The next section considers the sequential setting, in which the learned geometry and these regularity bounds are kept fixed while newly observed responses are incorporated into the repositories.
\section{Online continuation under fixed Lipschitz bounds}
\label{sec:update}

The global--local construction of \Cref{sec:global-local-operators} defines a predictor from fixed latent repositories. In a sequential setting, however, new responses become available over time and must be incorporated before subsequent predictions are made. Re-estimating the complete latent fit after every new observation would modify previously fitted values and potentially alter the regularity structure established during the structural fit. We instead keep the learned trajectory geometry and Lipschitz bounds fixed and update only the latent values associated with newly observed responses. This section formalises this predict--then--update mechanism, establishes that the fixed-bound continuation remains feasible, and compares it with exact insertion of the observed responses.

\subsection{Sequential repository update}
\label{sec:sequential-update}

At each time index $t$, all eligible groups are first predicted using only information available before that time. The global and local repositories remain unchanged while these predictions are computed, so no response observed at time $t$ can affect another prediction from the same time point. Only after the complete set of time-$t$ predictions has been stored are the corresponding responses revealed and used to update the repositories. The updated repositories are then used from time $t+1$ onwards.

Let $\mathcal G_t\subseteq \mathcal G$ denote the groups for which an eligible trajectory can be constructed at time $t$ and whose response is subsequently observed.
Before these observations are incorporated, the global latent repository can be written as
$
\mathcal R^G_{t-1}
=
\{(x_i,\widehat z_i^G,g_i):t_i<t\},
$
where $t_i$ denotes the target time associated with trajectory $x_i$.

Unlike the structural fit of \Cref{sec:robust-fit}, the online continuation does not re-estimate the geometry or the Lipschitz bounds. Its purpose is only to incorporate newly revealed responses while modifying them as little as necessary to preserve the fixed regularity structure. We therefore use an $L^1$ fidelity objective for the update.

The newly observed responses are not appended directly, since they may not satisfy the Lipschitz bound $\widehat\ell_G$ with the values already stored in the repository. Instead, we introduce a latent value $z^G_{g,t}$ for each new trajectory $x_{g,t}$.
For all groups updated at time $t$, we write $z_t^G=(z^G_{g,t})_{g\in\mathcal G_t}$.
The resulting latent values $\widehat z_t^G=(\widehat z^G_{g,t})_{g\in\mathcal G_t}$ are obtained jointly by solving
\begin{equation}
\label{eq:online-global-update}
\begin{aligned}
\widehat z_t^G
\in
\operatorname*{arg\,min}_{z_t^G\in\mathbb R^{|\mathcal G_t|}}
\quad&
\sum_{g\in\mathcal G_t}
\left|z^G_{g,t}-y_{g,t}\right|,
\\
\text{subject to}\quad&
\left|z^G_{g,t}-\widehat z_i^G\right|
\leq
\widehat\ell_G\,d_\omega(x_{g,t},x_i),
&&
g\in\mathcal G_t,\;
(x_i,\widehat z_i^G,g_i)\in\mathcal R^G_{t-1},
\\
&
\left|z^G_{g,t}-z^G_{h,t}\right|
\leq
\widehat\ell_G\,d_\omega(x_{g,t},x_{h,t}),
&&
g,h\in\mathcal G_t,\; g\neq h.
\end{aligned}
\end{equation}
The first set of constraints makes each new latent value compatible with all values already stored in the global repository. The second imposes the same Lipschitz condition between trajectories entering the repository at the same time point.

Once the problem has been solved, the global repository is updated as
$$
\mathcal R^G_t = \mathcal R^G_{t-1} \cup \left\{(x_{g,t},\widehat z^G_{g,t},g):g\in\mathcal G_t\right\},
$$
while previously stored latent values remain unchanged.

The $L^1$ objective does not in general guarantee uniqueness of the global batch update in \eqref{eq:online-global-update}. Using the same two-point data as in Remark~\ref{rem:l1-nonuniqueness}, let $y_1=0$ and $y_2=1$, assume that the previous repository is empty, and impose the fixed constraint $|z_1-z_2|\leq 1/2$. Then every pair $(z_1,z_2)=(a,a+1/2)$ with $a\in[0,1/2]$ minimises $|z_1|+|1-z_2|$. Thus, no uniqueness is assumed for the global batch update.

For the local repository, the pair structure defined in \Cref{sec:global-local-operators} contains only within-group constraints. Consequently, the online update separates into independent problems, one for each group. We define $z_t^L=(z^L_{g,t})_{g\in\mathcal G_t}$ and $\widehat z_t^L=(\widehat z^L_{g,t})_{g\in\mathcal G_t}$.
The new local latent values are obtained from
\begin{equation*}
\begin{aligned}
\widehat z_t^L
\in
\operatorname*{arg\,min}_{z_t^L\in\mathbb R^{|\mathcal G_t|}}
\quad&
\sum_{g\in\mathcal G_t}
\left|z^L_{g,t}-y_{g,t}\right|,
\\
\text{subject to}\quad&
\left|z^L_{g,t}-\widehat z_i^L\right|
\leq
\widehat\ell_L\,d_\omega(x_{g,t},x_i),
&&
g\in\mathcal G_t,\;
(x_i,\widehat z_i^L,g_i)\in\mathcal R^L_{t-1},
\;
g_i=g.
\end{aligned}
\end{equation*}
Unlike the global update, no constraints are imposed between new trajectories from different groups, so this problem separates into independent updates for each group. Once the new latent values have been obtained, they are appended to the local repository $\mathcal R^L_t$ in the same way as in the global case, while previously stored latent values remain unchanged.

\subsection{Feasibility, projection, and exact insertion}
\label{sec:online-properties}

\begin{proposition}%[Feasibility preservation under fixed Lipschitz bounds]
\label{prop:online-feasibility}
Suppose that the geometry $d_\omega$ and the global Lipschitz bound $\widehat\ell_G$ are fixed, and that the global repository $\mathcal R^G_{t-1}$ satisfies the Lipschitz condition with bound $\widehat\ell_G$. Then the feasible set of Problem~\ref{eq:online-global-update} is nonempty and the minimum is attained. Moreover, appending any feasible solution for the new latent values preserves the same Lipschitz bound $\widehat\ell_G$ on the enlarged repository $\mathcal R^G_t$.

The same conclusion holds for each group-specific local repository under the fixed shared bound $\widehat\ell_L$, applying the argument separately within each group.
\end{proposition}
\begin{proof}
Consider first the global repository. If $\mathcal R^G_{t-1}$ is empty, the old--new constraints are absent and feasibility is immediate by assigning the same constant value to all new latent variables, which also satisfies every new--new Lipschitz constraint. Otherwise, since $\mathcal R^G_{t-1}$ satisfies the Lipschitz condition with bound $\widehat\ell_G$, the McShane--Whitney theorem provides an $\widehat\ell_G$-Lipschitz extension of its latent values to the whole trajectory space. Evaluating any such extension at all new trajectories $x_{g,t}$, $g\in\mathcal G_t$, provides a simultaneous assignment of new latent values satisfying both the old--new constraints and the new--new constraints in Problem~\ref{eq:online-global-update}. Hence the feasible set is nonempty.

The global update minimises
$$
\sum_{g\in\mathcal G_t}
\left|z^G_{g,t}-y_{g,t}\right|
$$
over this nonempty feasible set. As in the proof of Proposition~\ref{prop:latent-fit}, continuity and coercivity of the objective, together with closedness of the feasible set, imply that the minimum is attained.

The old repository already satisfies all old--old Lipschitz inequalities with bound $\widehat\ell_G$. Any feasible solution of the update satisfies all required old--new and new--new inequalities with the same bound. Therefore, after appending the new latent values, the enlarged repository $\mathcal R^G_t$ also satisfies the Lipschitz condition with bound $\widehat\ell_G$.

For the local repository, the same argument is applied independently within each group using the fixed bound $\widehat\ell_L$.
\end{proof}

\begin{remark}
For a single new trajectory $x$ with observed response $y$ and a nonempty previous repository satisfying the fixed Lipschitz bound $L$, the feasible set of the online update is exactly the McShane--Whitney interval $[f_L^M(x),f_L^W(x)]$. Under the $L^1$ objective, the new latent value is therefore the projection of $y$ onto this interval. If the relevant previous repository is empty, there are no old--new constraints and the update gives $\widehat z=y$.
\end{remark}

The fixed-bound update can also be contrasted with exact online interpolation, in which every newly observed response is inserted without modification and the empirical Lipschitz constant is allowed to adapt to the enlarged repository.
Let $\mathcal R_{t-1}^{\mathrm{raw}}$ denote the exact repository available before time $t$, and let $\mathcal N_t^{\mathrm{raw}}= \{(x_{g,t},y_{g,t},g):g\in\mathcal G_t\}$ denote the set of newly observed trajectories and responses inserted at time $t$. The enlarged exact repository is therefore
$$
\mathcal R_t^{\mathrm{raw}}
=
\mathcal R_{t-1}^{\mathrm{raw}}
\cup
\mathcal N_t^{\mathrm{raw}}.
$$
For any required pair $(i,j)$, its raw pairwise slope is defined as
$$
r_{ij}^{\mathrm{raw}}
=
\begin{cases}
\displaystyle
\frac{|y_i-y_j|}{d_\omega(x_i,x_j)},
& d_\omega(x_i,x_j)>0,\\[3mm]
0,
& d_\omega(x_i,x_j)=0 \text{ and } y_i=y_j,\\[1mm]
+\infty,
& d_\omega(x_i,x_j)=0 \text{ and } y_i\neq y_j.
\end{cases}
$$

Under a fixed trajectory geometry and a fixed rule determining which pairs are Lipschitz-constrained, let $L_{t-1}^{\mathrm{raw}}$ be the maximum raw slope over the required pairs contained in $\mathcal R_{t-1}^{\mathrm{raw}}$, let $L_{\mathrm{old-new},t}^{\mathrm{raw}}$ be the maximum over the required pairs connecting $\mathcal R_{t-1}^{\mathrm{raw}}$ with $\mathcal N_t^{\mathrm{raw}}$, and let $L_{\mathrm{new-new},t}^{\mathrm{raw}}$ be the maximum over the required pairs contained in $\mathcal N_t^{\mathrm{raw}}$. Let $L_t^{\mathrm{raw}}$ denote the maximum raw slope over all required pairs in the enlarged repository $\mathcal R_t^{\mathrm{raw}}$, with the maximum of an empty family defined as zero. Since every required pair in $\mathcal R_t^{\mathrm{raw}}$ belongs to exactly one of the old--old, old--new, and new--new classes,
\begin{equation*}
L_t^{\mathrm{raw}}
=
\max\left\{
L_{t-1}^{\mathrm{raw}},
L_{\mathrm{old-new},t}^{\mathrm{raw}},
L_{\mathrm{new-new},t}^{\mathrm{raw}}
\right\}.
\end{equation*}
In particular, $L_t^{\mathrm{raw}}\geq L_{t-1}^{\mathrm{raw}}$, so along nested exact repositories with fixed geometry and a fixed pair rule, the exact empirical Lipschitz constant is nondecreasing. It may remain unchanged when the newly added observations introduce no steeper required pair, but it cannot decrease.

Returning to the projected online update, this procedure is different from refitting the structural latent model after each new time point. The previously stored latent values remain unchanged, the trajectory geometry and the bounds $\widehat\ell_G$ and $\widehat\ell_L$ remain fixed, and only the latent values associated with the newly observed responses are estimated through the $L^1$ objective before being appended to the repositories. Rather than re-estimating the Lipschitz bounds as new observations become available, the new latent values are chosen so that the bounds obtained from the structural fit remain valid after the repositories are enlarged.

This predict--then--update cycle can be iterated for an arbitrary sequence of observation times. Starting from feasible global and local latent repositories, \Cref{prop:online-feasibility} guarantees inductively that the fixed Lipschitz bounds remain valid after every update, while previously fitted latent values and structural parameters remain unchanged. The resulting procedure therefore provides a fixed-bound sequential continuation of the fitted Lipschitz model rather than a repeated refitting of that model. The computational realisation of the structural fit and online updates is described in the next section.
\section{Computational implementation}
\label{sec:computational-implementation}

The mathematical formulation developed in the previous sections leads to finite-dimensional convex optimisation problems, but their direct solution may become expensive as the number of trajectories increases. Two computational observations are particularly useful. First, during the sequential continuation of \Cref{sec:update}, the trajectory geometry and all previously stored trajectories remain fixed. Distances between previously stored trajectories therefore do not need to be recomputed; at each update, only distances involving the newly added trajectories are required. Second, the structural Lipschitz fit involves a potentially quadratic number of pairwise constraints, which motivates the constraint-generation procedure described below.

For a fixed trajectory geometry, the Lipschitz fitting problem requires $|z_i-z_j|\leq \ell\,d_\omega(x_i,x_j)$ for every admissible pair of trajectories. As the number $n$ of available trajectories increases, the number of pairwise constraints can grow as $\mathcal{O}(n^2)$, making the explicit construction and simultaneous storage of the complete constraint system increasingly expensive.

We therefore use a constraint-generation procedure, following the general principle of solving successively enlarged reduced problems rather than imposing all pairwise constraints from the outset \citep{lin2022constraint}. The optimisation is first solved with only a subset $\mathcal A_0$ of the admissible pairwise constraints. Let $(\widehat z^{(m)},\widehat\ell^{(m)})$ denote the solution obtained at iteration $m$. The remaining pairs are then checked against the Lipschitz condition, and those violating it beyond a numerical tolerance $\varepsilon\geq0$ are collected in
$$
\mathcal V_m^{\varepsilon}
=
\left\{
(i,j)\in\mathcal P\setminus\mathcal A_m:
\left|
\widehat z_i^{(m)}-\widehat z_j^{(m)}
\right|
>
\widehat\ell^{(m)}d_\omega(x_i,x_j)+\varepsilon
\right\}.
$$
Violated constraints are added to the current system, so that $\mathcal A_{m+1}=\mathcal A_m\cup \mathcal V_m^{\varepsilon}$, and the optimisation is solved again. The process continues until no admissible pair violates the Lipschitz condition beyond the prescribed tolerance; that is, until the first iteration $M$ for which $\mathcal V_M^\varepsilon=\emptyset$.

The preceding procedure has a simple exact guarantee when the reduced problems are solved exactly and constraint separation is performed without numerical tolerance. The following result records this idealised property.

\begin{proposition}%[Correctness and finite termination of exact constraint generation]
\label{prop:constraint-generation}
Let $\mathcal P$ be a finite family of admissible pairs and let $\mathcal A_0\subseteq\mathcal P$. At iteration $m$, suppose that Problem~\ref{pr:latent-fit} is solved exactly with the constraints restricted to $\mathcal A_m$, yielding a minimiser $(\widehat z^{(m)},\widehat\ell^{(m)})$. Define
$$
\mathcal V_m
=
\left\{
(i,j)\in\mathcal P\setminus\mathcal A_m:
\left|
\widehat z_i^{(m)}-\widehat z_j^{(m)}
\right|
>
\widehat\ell^{(m)}d_\omega(x_i,x_j)
\right\}.
$$
If $\mathcal V_m=\varnothing$, stop. Otherwise set $\mathcal A_{m+1}=\mathcal A_m\cup\mathcal V_m$ and repeat. Then the procedure terminates after finitely many iterations. At termination, the resulting pair $(\widehat z^{(M)},\widehat\ell^{(M)})$ is a minimiser of the full problem with constraint family $\mathcal P$.
\end{proposition}
\begin{proof}
Whenever $\mathcal V_m\neq\varnothing$, at least one previously inactive constraint is added, so $\mathcal A_m\subsetneq\mathcal A_{m+1}$. Since $\mathcal P$ is finite, this can occur only finitely many times, and the procedure must terminate at some iteration $M$ with $\mathcal V_M=\varnothing$. At that point, $(\widehat z^{(M)},\widehat\ell^{(M)})$ satisfies every constraint in $\mathcal P$ and is therefore feasible for the full problem. On the other hand, because $\mathcal A_M\subseteq\mathcal P$, the feasible set of the restricted problem contains the feasible set of the full problem, so the optimal value of the restricted problem cannot exceed the optimal value of the full problem. Since the terminal restricted minimiser is itself feasible for the full problem, its objective value cannot be smaller than the full optimum. The two objective values are therefore equal, and $(\widehat z^{(M)},\widehat\ell^{(M)})$ is a minimiser of the full problem.
\end{proof}

\begin{remark}
The same conclusion holds if only a nonempty subset of the currently violated constraints is added at each nonterminal iteration, provided that exact separation is performed over the full finite family $\mathcal P$. In that case, at least one new constraint is still added whenever the current solution is not feasible for the full problem.
\end{remark}

Proposition~\ref{prop:constraint-generation} concerns an exact scheme. In the numerical implementation, the reduced convex problems are solved in finite precision, and a tolerance $\varepsilon>0$ is used to assess constraint violations. After termination, the computed solution is checked against the complete family of admissible pairs, including the compatibility conditions associated with zero-distance pairs, and accepted only when the maximum residual violation is within the prescribed tolerance. Thus, constraint generation remains a computational device that does not modify the underlying Lipschitz-constrained problem, and Proposition~\ref{prop:constraint-generation} should not be interpreted as an exact finite-precision convergence guarantee.

\section{Application to sow-production trajectories}
\label{sec:application}

\begin{figure}[htbp]
\centering
\resizebox{\linewidth}{!}{%
\begin{tikzpicture}[x=1cm,y=1cm,>=stealth]
  % Local colours and styles: no definitions elsewhere in the article are needed.
  \definecolor{lfnavy}{RGB}{29,55,87}
  \definecolor{lfblue}{RGB}{37,132,215}
  \definecolor{lforange}{RGB}{234,145,33}
  \definecolor{lfgreen}{RGB}{62,170,100}
  \definecolor{lfpurple}{RGB}{133,87,193}
  \definecolor{lfgray}{RGB}{170,188,207}
  \definecolor{lfred}{RGB}{216,59,91}
  \definecolor{lfviolet}{RGB}{105,55,175}
  \tikzset{
    lfpanel/.style={rounded corners=5pt,line width=.75pt},
    lfarrow/.style={->,draw=lfnavy,line width=1.05pt,shorten >=1pt},
    lfheading/.style={font=\sffamily\bfseries\fontsize{8}{9.5}\selectfont,align=center},
    lfmath/.style={font=\fontsize{8.0}{9.5}\selectfont,align=center},
    lfmini/.style={font=\fontsize{7}{8.5}\selectfont,align=center},
  }

  % Five main panels: total drawing width is 14.90 cm.
  \draw[lfpanel,draw=lfgray,fill=lfgray!9] (0,0) rectangle (1.90,3.90);
  \draw[lfpanel,draw=lfred!60,fill=lfred!5] (2.20,0) rectangle (4.38,3.90);
  \draw[lfpanel,draw=lfblue!70,fill=lfblue!6] (4.68,0) rectangle (8.30,3.90);
  \draw[lfpanel,draw=lfgreen!70,fill=lfgreen!6] (8.60,0) rectangle (12.22,3.90);
  \draw[lfpanel,draw=lforange!80,fill=lforange!7] (12.52,0) rectangle (14.90,3.90);

  % Panel headings (all contained inside their panels).
  \node[lfheading] at (.95,3.5) {Trajectory\\data};
  \node[lfheading] at (3.29,3.35) {Weighted\\trajectory\\geometry};
  \node[lfheading] at (6.49,3.5) {Robust Lipschitz\\regression};
  \node[lfheading] at (10.41,3.59) {Global--local prediction};
  \node[lfheading] at (13.71,3.5) {Validation\\error};

  % ---- 1: trajectory data (three schematic series) ----
    \begin{scope}[shift={(0.15,0.9)}]
      \draw[lfblue,line width=1.2pt] plot[smooth] coordinates
         {(.13,1.48)(.39,1.69)(.71,2.02)(1.03,1.64)(1.35,1.72)(1.58,1.90)};
      \foreach \a/\b in {.13/1.48,.39/1.69,.71/2.02,1.03/1.64,1.35/1.72,1.58/1.90}
        \fill[lfblue] (\a,\b) circle (2.4pt);
    
      \draw[lforange,line width=1.2pt] plot[smooth] coordinates
         {(.13,.88)(.39,1.07)(.71,1.17)(1.03,.95)(1.35,.81)(1.58,1.05)};
      \foreach \a/\b in {.13/.88,.39/1.07,.71/1.17,1.03/.95,1.35/.81,1.58/1.05}
        \fill[lforange] (\a,\b) circle (2.4pt);
    
      \draw[lfgreen,line width=1.2pt] plot[smooth] coordinates
         {(.13,.16)(.39,.36)(.71,.50)(1.03,.31)(1.35,.16)(1.58,.38)};
      \foreach \a/\b in {.13/.16,.39/.36,.71/.50,1.03/.31,1.35/.16,1.58/.38}
        \fill[lfgreen] (\a,\b) circle (2.4pt);
    
      \node[lfmini] at (.89,-.12) {$\cdots$};
    
      % time axis
      \draw[->,draw=lfnavy,line width=.7pt] (0.2,-0.34)--(1.52,-0.34);
      \node[lfmini] at (.86,-0.56) {Time};
    \end{scope}

  % ---- 2: learned weighted geometry ----
    \node[lfmath] at (3.29,0.31) {$d_{\omega}(x_i,x_j)$};
    
    \begin{scope}[shift={(2.30,0.56)}]
      % named points
      \coordinate (A)  at (0.26,1.96);   % purple
      \coordinate (B)  at (1.06,2.05);   % top light blue
      \coordinate (C)  at (1.68,1.85);   % top-right blue
      \coordinate (D)  at (0.12,1.2);   % left blue
      \coordinate (E)  at (0.65,0.62);   % lower-left blue
      \coordinate (F)  at (1.16,0.20);   % green
      \coordinate (G)  at (1.77,1.30);   % orange
      \coordinate (H)  at (1.76,0.60);   % orange
      \coordinate (M1) at (0.82,1.5);   % gray
      \coordinate (M2) at (1.28,1.3);   % gray
      \coordinate (M3) at (0.98,0.9);   % gray
      % Dashed weighted-geometry links
        \foreach \u/\v in {A/B,A/M1,A/M2,B/M2,D/M1,E/M1,E/M3,M1/M2,M1/M3,M2/M3,M2/G,M3/G,M2/H,G/H,F/M3}{
            \draw[lfgray!75,dashed,line width=.8pt] (\u)--(\v);
        }
      % nodes
      \fill[lfpurple] (A) circle (3.4pt);
      \fill[lfblue!55] (B) circle (3.2pt);
      \fill[lfblue!65] (C) circle (3.2pt);
      \fill[lfblue] (D) circle (3.2pt);
      \fill[lfblue] (E) circle (3.2pt);
      \fill[lfgreen] (F) circle (3.2pt);
      \fill[lforange] (G) circle (3.2pt);
      \fill[lforange] (H) circle (3.2pt);
      \fill[lfgray] (M1) circle (2.8pt);
      \fill[lfgray] (M2) circle (2.8pt);
      \fill[lfgray] (M3) circle (2.8pt);
    \end{scope}

  % ---- 3: robust latent fit and a stylised smoother ----
    \node[lfmini] at (6.49,2.35)
      {$\displaystyle \min_z\!\left\{\!\begin{gathered}
       \frac1n\sum_i\rho(y_i-z_i)\\[-1.5pt]
       +\eta\,\mathrm{Lip}_{\mathcal P}(z)
       \end{gathered}\!\right\}$};
    
    \begin{scope}[shift={(4.96,.18)}]
      \draw[lfblue!16,line width=6pt,line cap=round]
        plot[smooth] coordinates {(.18,.28)(.66,.60)(1.23,.96)(1.82,.78)(2.45,.56)(2.95,.74)};
      \draw[lfblue,line width=1.4pt]
        plot[smooth] coordinates {(.18,.28)(.66,.60)(1.23,.96)(1.82,.78)(2.45,.56)(2.95,.74)};
    
      \foreach \x/\y in {
        .26/.14,.50/.82,.76/.72,1.00/1.08,1.24/.60,
        1.52/1.14,1.78/.48,2.10/.90,2.42/.40,2.72/.80
      }{
        \fill[lfgray] (\x,\y) circle (2.9pt);
      }
    \end{scope}

  % ---- 4: global--local block and final mixture (single combined panel) ----
    \node[lfmini] at (10.41,.9)
      {$\widehat Y(x)=(1-\alpha)\widehat Y_G(x)$};
    \node[lfmini] at (10.41,.36)
      {$+\alpha\widehat Y_L(x)$};
    
    % Local component: left
    \begin{scope}[shift={(8.86,2.25)}]
      \draw[dashed,draw=lfgray,line width=.65pt] (.42,.58) circle (.28);
      \foreach \x/\y in {.20/.68,.36/.86,.56/.77,.21/.42,.49/.43,.65/.57}
        \fill[lforange!85] (\x,\y) circle (2.75pt);
      \draw[->,draw=lfnavy!75,line width=.8pt] (.82,.58)--(1.15,.58);
      \fill[lfred] (1.22,.58) circle (3.3pt);
    \end{scope}
    
    % Global component: right
    \begin{scope}[shift={(8.86,1.44)}]
      \foreach \x/\y in {.12/.26,.34/.58,.60/.32,.80/.64}
        \fill[lfblue!55] (\x,\y) circle (2.5pt);
      \foreach \x/\y in {.94/.33,1.14/.68,1.33/.38}
        \fill[lfgreen!65] (\x,\y) circle (2.5pt);
      \foreach \x/\y in {.21/.66,.49/.37,1.28/.53}
        \fill[lfgray] (\x,\y) circle (2.5pt);
    
      \draw[lfblue,line width=1pt]
        plot[smooth] coordinates {(.02,.38)(.46,.48)(.96,.38)(1.42,.57)};
      \fill[lfviolet] (1.43,.56) circle (3.3pt);
    \end{scope}
    
    \node[lfmath] at (11.39,2.82) {$\widehat Y_L(x)$};
    \node[lfmath] at (11.39,1.94) {$\widehat Y_G(x)$};

  % ---- 5: validation error ----
  \node[lfmath] at (13.71,2.74) {$E_{\mathrm{sel}}(\omega)$};
    \begin{scope}[shift={(12.78,.72)}]
      \draw[->,draw=lfnavy,line width=.7pt] (0,0)--(1.84,0);
      \draw[->,draw=lfnavy,line width=.7pt] (0,0)--(0,1.58);
    
      \filldraw[fill=lfgray!60,draw=lfgray] (.26,0) rectangle (.53,.72);
      \filldraw[fill=lfgray!80,draw=lfgray] (.72,0) rectangle (.99,1.08);
      \filldraw[fill=lfgray!60,draw=lfgray] (1.18,0) rectangle (1.45,.58);

      \node[lfmini] at (.92,-.35) {Validation folds};
    \end{scope}

  % ---- Forward pipeline and outer optimisation return ----
  \foreach \a/\b in {1.90/2.20,4.38/4.68,8.30/8.60,12.22/12.52}
    \draw[lfarrow] (\a,1.92)--(\b,1.92);

  \draw[lfpanel,draw=lfnavy!70,fill=lfgray!13]
       (5.1,4.43) rectangle (11.89,6.0);
  \node[lfheading] at (7.56,5.62) {Outer $\omega$-optimisation};
  \node[lfmath,text width=4.1cm] at (7.1,4.95)
     {$\displaystyle
       \widehat\omega\approx
       \operatorname*{arg\,min}_{\omega\in\Delta_p}
       E_{\mathrm{sel}}(\omega)$};

  % Illustrative one-dimensional slice through candidate errors.
    \begin{scope}[shift={(9.90,4.67)}]
      \draw[->,draw=lfnavy,line width=.60pt] (0,0)--(1.36,0);
      \draw[->,draw=lfnavy,line width=.60pt] (0,0)--(0,1.00);
    
      \node[lfmini] at (1.52,0) {$\omega$};
      \node[lfmini] at (.54,0.97) {$E_{\mathrm{sel}}(\omega)$};
    
      \draw[lfblue!45,line width=1.15pt]
        plot[smooth] coordinates
          {(.08,.78)(.22,.56)(.40,.34)(.63,.22)(.88,.28)(1.10,.50)(1.26,.78)};
    
      \foreach \x/\y in {
        .08/.78,.16/.66,.22/.56,.31/.44,.40/.34,.51/.26,
        .63/.22,.77/.24,.88/.28,1.00/.38,1.10/.50,1.18/.64,1.26/.78
      }{
        \fill[lfblue!55] (\x,\y) circle (1.6pt);
      }
    
      \fill[lfpurple] (.63,.22) circle (2.5pt);
    \end{scope}

  % Outer search sends new candidate weights to geometry only.
  \draw[lfarrow] (5.1,5.215)--(3.28,5.215)--(3.28,3.90);
  % Validation errors feed back to outer search (not to input data).
  \draw[lfarrow] (13.71,3.90)--(13.71,5.215)--(11.89,5.215);
\end{tikzpicture}
}
\caption{Schematic representation of the proposed regression framework.}
\label{fig:regression-framework}
\end{figure}
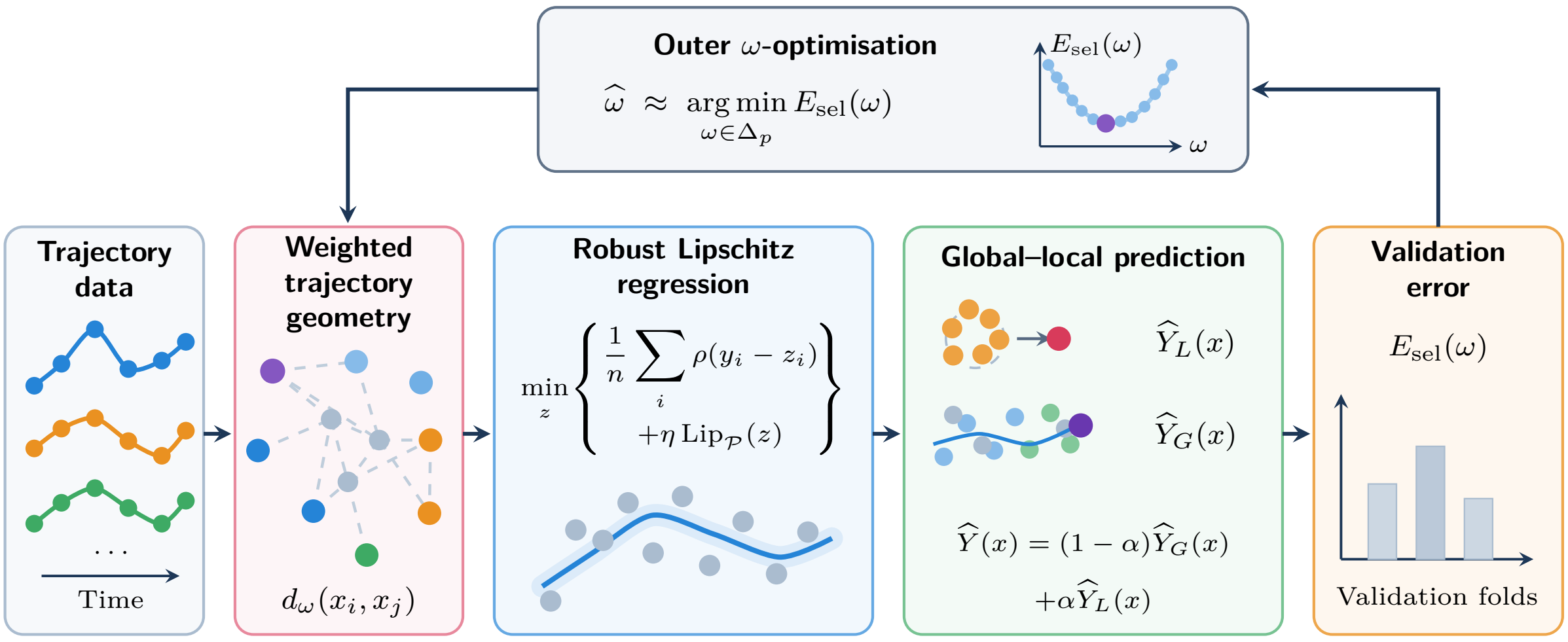

The purpose of this section is to provide a compact empirical illustration of the mathematical framework developed above in a real longitudinal prediction problem. Rather than performing a comprehensive applied study, we use a single sow-production outcome to instantiate the trajectory pseudometric, the grouped global--local construction, and the fixed-bound sequential continuation, and to illustrate some of the geometric effects established in the preceding sections. A broader applied analysis is beyond the scope of the present mathematical study and is reserved for a separate applied study.
The overall predictive architecture is illustrated in \Cref{fig:regression-framework}.

\subsection{Problem setup and temporal evaluation protocol}

We consider a longitudinal forecasting problem based on monthly climatic and production records from $52$ commercial pig farms in Spain, spanning approximately ten years, from 2014 to 2023. The objective is to predict the monthly farrowing rate of each farm, represented by \texttt{FarrowingRate119}, from its recent climatic and productive history. Specifically, the predictors are mean temperature, mean relative humidity, and previous values of the farrowing rate itself. Since sow gestation typically lasts approximately $114$--$116$ days, close to four months, we represent the recent history of each farm through the four months preceding the prediction month. This temporal window is therefore consistent with the biological time scale underlying the reproductive outcome considered here.

This problem fits naturally into the trajectory-based framework developed in the preceding sections. Let $\mathcal F$ denote the set of farms. In this application, each farm constitutes one group, so that $\mathcal G=\mathcal F$. For farm $f\in\mathcal F$ and month $t$, let $T_{f,t}$ denote mean temperature, $H_{f,t}$ mean relative humidity, and $Y_{f,t}$ the farrowing-rate response. The monthly state is $s_{f,t}=(T_{f,t},H_{f,t},Y_{f,t})\in\mathbb R^3$, and the trajectory associated with the prediction of $Y_{f,t}$ is formed by the four preceding states,
$$
x_{f,t}
=
\left(
s_{f,t-4},
s_{f,t-3},
s_{f,t-2},
s_{f,t-1}
\right)
\in\mathbb R^{12}.
$$
Thus, in the notation of \Cref{sec:trajectory-spaces}, $q=3$, $k=4$, and $p=qk=12$. Before computing trajectory distances, the state variables are standardised so that the coordinates entering $d_\omega$ are on a common scale. Farm identity is not a coordinate of $x_{f,t}$ and therefore does not enter the learned pseudometric $d_\omega$; instead, it plays the role of the group label $g=f$, determining the global and local pair structures introduced in \Cref{sec:global-local-structures}.

Model development is restricted to the pre-2021 portion of the data, corresponding approximately to the first seven years of the longitudinal panel. The final model specification---including the trajectory weights $\omega$, the Lipschitz regularisation parameter $\eta$, the McShane--Whitney mixing coefficients $\lambda_G$ and $\lambda_L$, and the global--local coefficient $\alpha$---is selected using an expanding-window form of time-series cross-validation with four annual folds (2017--2020). In each fold, the model is fitted using only observations preceding the validation year and is then evaluated sequentially throughout that year. The training period therefore grows from one fold to the next, while validation always uses subsequent observations.

The trajectory weights $\omega$ and regularisation parameter $\eta$ are searched using the Tree-structured Parzen Estimator (TPE) implemented in Optuna \citep{akiba2019optuna}, while the prediction-mixing coefficients are calibrated by grid search. In this application, the selection error $E_{\mathrm{sel}}$ shown in \Cref{fig:regression-framework} is the average of the farm-macro mean absolute errors over the four validation folds. The configuration with the lowest selection error is retained.

Once selected, the model and its structural parameters are fixed and carried forward to the 2021--2023 period, where predictions are generated sequentially using the predict--then--update procedure described in \Cref{sec:sequential-update}, without further parameter selection. Within each evaluation month, all eligible farms are predicted before any response from that month is incorporated; the subsequent global update is performed jointly across farms, whereas the local update separates by farm.

\subsection{Predictive and geometric results}

The final Lipschitz predictor selected by the pre-2021 validation protocol is summarised in Table~\ref{tab:selected-lipschitz-configuration}. In particular, it uses the robust pair-restricted Lipschitz regression of \Cref{sec:robust-fit} with Huber fidelity for the latent-response fit, together with the global and local McShane--Whitney prediction branches of \Cref{sec:global-local-operators}.

\begin{table}[!htbp]
\centering
\caption{Configuration of the selected Lipschitz predictor.}
\label{tab:selected-lipschitz-configuration}
\begin{tabular}{ll}
\toprule
Component & Selected specification \\
\midrule
State variables & $T$, $H$, $Y$(\textit{farrowing rate}) \\
Trajectory length & $k=4$ months \\
Trajectory dimension & $p=12$ \\
Latent-response fidelity & Huber with $\delta=1$ \\
Trajectory weights & selected by pre-2021 validation \\
Lipschitz regularisation & $\eta=0.075$ \\
Global McShane--Whitney coefficient & $\lambda_G=0.515625$ \\
Local McShane--Whitney coefficient & $\lambda_L=0.484375$ \\
Global--local coefficient & $\alpha=0.2625$ \\
\bottomrule
\end{tabular}
\end{table}

To place its predictive performance in context, we compare this selected Lipschitz predictor with three standard supervised-learning benchmarks and with two simpler instances of the Lipschitz construction. Random Forests and neural networks are widely used nonlinear approaches in data-driven livestock prediction problems \citep{sanjuan2025soft}, and are included here as global reference models. Both are fitted using observations from all farms and receive the same $12$ standardised trajectory coordinates as the proposed predictor, together with a one-hot encoding of farm identity, allowing farm-specific differences to be represented within a single model shared across farms. We additionally consider a local linear regression benchmark, for which a separate ordinary least-squares model is fitted to each farm using only that farm's previous observations and the same trajectory coordinates.

The two additional Lipschitz configurations are included to illustrate the effect of the main methodological steps developed in this article. The first is an exact local Lipschitz predictor, which uses the learned trajectory geometry but directly interpolates the observed responses, without introducing latent values or a regularisation penalty. The second retains the local information scope but replaces exact interpolation by the robust latent formulation of \Cref{sec:robust-fit}, using Huber fidelity. The selected model then additionally incorporates the global--local construction of \Cref{sec:global-local-operators}. All models are evaluated over the same 2021--2023 period under the same chronological information constraint and month-by-month prediction protocol. Table~\ref{tab:predictive-comparison} reports the resulting farm-macro MAE.

\begin{table}[!htbp]
\centering
\caption{Farm-macro MAE over the final 2021--2023 evaluation period.}
\label{tab:predictive-comparison}
\begin{tabular}{lc}
\toprule
Model & Farm-macro MAE \\
\midrule
Selected robust Lipschitz model & $3.2973$ \\
Global neural network & $3.3603$ \\
Global Random Forest & $3.3662$ \\
Robust local Lipschitz (Huber) & $3.4164$ \\
Local linear regression & $3.5009$ \\
Exact local Lipschitz & $3.8971$ \\
\bottomrule
\end{tabular}
\end{table}

The comparison illustrates the predictive changes associated with the main modelling steps developed in this work on the final evaluation period. Moving from exact local interpolation to the local Huber latent formulation reduces the farm-macro MAE from $3.8971$ to $3.4164$, while the final global--local construction further reduces it to $3.2973$. The selected robust Lipschitz predictor therefore attains the lowest observed MAE among the models reported in Table~\ref{tab:predictive-comparison}, including the global Random Forest and neural-network references. These comparisons are intended as an empirical illustration of the framework rather than as evidence of general predictive superiority.

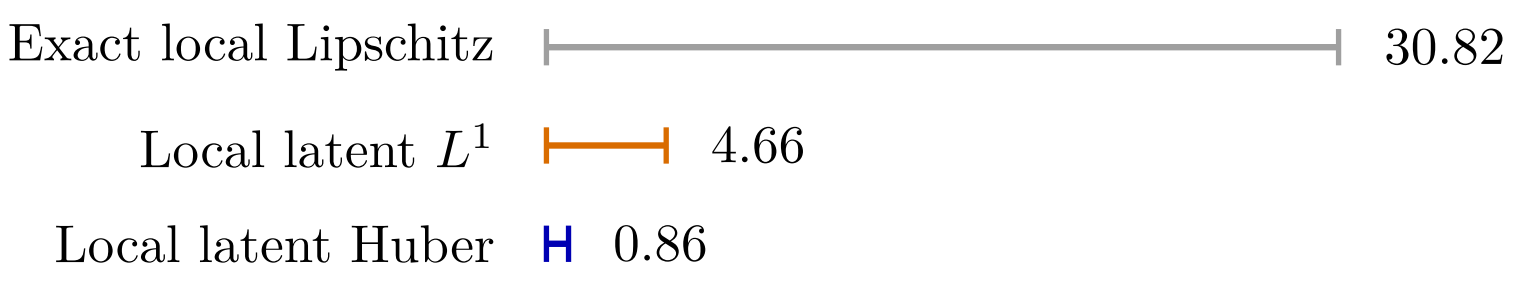
\begin{figure}[htbp]
\centering
\begin{tikzpicture}[x=0.17cm,y=1cm]
% Exact local Lipschitz
\node[anchor=east,font=\small] at (-1.2,1.8) {Exact local Lipschitz};
\draw[gray!75,line width=1.2pt] (0,1.8) -- (30.82,1.8);
\draw[gray!75,line width=1pt] (0,1.68) -- (0,1.92);
\draw[gray!75,line width=1pt] (30.82,1.68) -- (30.82,1.92);
\node[anchor=west,font=\small] at (31.8,1.8) {$30.82$};
% Local latent L1
\node[anchor=east,font=\small] at (-1.2,1.15) {Local latent $L^1$};
\draw[orange!85!black,line width=1.2pt] (0,1.15) -- (4.66,1.15);
\draw[orange!85!black,line width=1pt] (0,1.03) -- (0,1.27);
\draw[orange!85!black,line width=1pt] (4.66,1.03) -- (4.66,1.27);
\node[anchor=west,font=\small] at (5.6,1.15) {$4.66$};
% Local latent Huber
\node[anchor=east,font=\small] at (-1.2,0.50) {Local latent Huber};
\draw[blue!70!black,line width=1.2pt] (0,0.50) -- (0.86,0.50);
\draw[blue!70!black,line width=1pt] (0,0.38) -- (0,0.62);
\draw[blue!70!black,line width=1pt] (0.86,0.38) -- (0.86,0.62);
\node[anchor=west,font=\small] at (1.8,0.50) {$0.86$};
\end{tikzpicture}
\caption{Mean local McShane--Whitney extension width for the three local predictors considered in the application.}
\label{fig:local-extension-width}
\end{figure}

A complementary geometric diagnostic is provided by the local McShane--Whitney extension width, summarised in \Cref{fig:local-extension-width}. For each pre-2021 prediction, we compute the width $f^{W}_{L}(x)-f^{M}_{L}(x)$ of the corresponding local extension interval and then average it over predictions.

The large width obtained under exact interpolation reflects the considerable freedom left by the corresponding Lipschitz constraints when extending the observed responses to a new trajectory. In particular, exact interpolation forces all observed responses to be reproduced, so isolated steep pairwise slopes---which may arise from measurement noise, atypical observations, unobserved farm conditions, or genuine production variability---can increase the empirical Lipschitz bound and consequently weaken the restriction imposed on the extension. The resulting McShane--Whitney interval can therefore be broad, leaving a wide range between the two extremal admissible extensions from which the prediction rule selects its value. Replacing exact interpolation by latent Lipschitz fitting is accompanied by a substantial contraction of this range, producing a considerably more restrictive extension geometry. This geometric reduction complements the predictive differences reported in Table~\ref{tab:predictive-comparison}.

\subsection{Illustration of the mathematical framework}

The sow-production application provides a concrete instantiation of the mathematical construction developed above. Each predictive input is represented by a trajectory in the learned pseudometric space $(\mathcal X,d_\omega)$, while farm identity determines the group structure without entering the trajectory coordinates themselves. The subsequent prediction therefore depends on both the learned notion of similarity between recent production histories and the global--local organisation of the available observations.

The extension-width results illustrate the role of Lipschitz regularisation. Exact interpolation can require a large Lipschitz bound when nearby trajectories have substantially different observed responses. Robust latent fitting relaxes this requirement by allowing deviations from the observations while penalising Lipschitz variation. In particular, \Cref{thm:huber-strict-shrinkage} establishes strict Lipschitz shrinkage under Huber fidelity and its stated assumptions. Since \Cref{eq:rel_width_L} gives $0\leq f^W_L(x)-f^M_L(x)\leq 2L\,d(x,A)$, smaller Lipschitz bounds provide tighter upper control on extension widths. This helps interpret the narrower intervals observed for the latent models, although it does not guarantee smaller realised widths or improved predictive accuracy.

The selected predictor further illustrates the grouped construction of \Cref{sec:global-local-operators}. The global and local McShane--Whitney operators use different anchor sets, representing respectively cross-farm and same-farm information, and their predictions are combined through the coefficient $\alpha$. The interior value selected in the application therefore provides a concrete example in which both branches contribute to the final predictor. Importantly, this combination remains mathematically controlled: the regularity of the resulting prediction is bounded in terms of the Lipschitz bounds of the two component branches.

Finally, the month-by-month evaluation protocol implements the fixed-bound sequential continuation described in \Cref{sec:sequential-update}. Once the trajectory geometry, regularisation parameters, extension coefficients, and Lipschitz bounds have been fixed, predictions are generated from the current repositories before the new responses are revealed, and only the newly observed trajectories are subsequently incorporated. The application therefore illustrates the complete sequence developed in this work: learned trajectory geometry, robust pair-restricted Lipschitz fitting, McShane--Whitney extension, global--local combination, and causal sequential continuation.
\section{Conclusions}
\label{sec:conclusions}

This work develops a robust Lipschitz regression framework for grouped longitudinal data represented as trajectories in a possibly learned pseudometric space. The construction combines a pair-restricted variational fit of latent responses, McShane--Whitney extension, global and local information structures, and a fixed-bound sequential continuation. Its main mathematical contribution lies in the structural properties that emerge from this combination. For a fixed trajectory geometry, positive Huber--Lipschitz regularisation strictly reduces the pair-restricted Lipschitz seminorm whenever the observed response has finite positive seminorm, whereas $L^1$ fidelity can preserve exact interpolation below a data-dependent threshold. In the grouped setting, the shared local Lipschitz bound induces an exact threshold-and-saturation mechanism across groups, while the global--local predictor retains an explicit regularity bound determined by those of its two component branches. The sequential continuation then preserves the fitted regularity structure as new observations arrive, and for a single new observation reduces to projection onto the corresponding McShane--Whitney interval.

The framework also provides a direct geometric interpretation of the different stages of prediction. The learned pseudometric determines which trajectories are considered similar, the fitted Lipschitz bounds control the admissible variation of the latent response over that geometry, and the McShane--Whitney operators determine the extremal extensions compatible with the fitted repositories. Group identity is deliberately kept outside the trajectory geometry and instead determines the global and local information scopes, while the coefficient $\alpha$ controls the algebraic balance between the resulting prediction branches. These quantities should be interpreted predictively rather than causally: the trajectory weights describe relevance within the selected geometry, $\alpha$ is not a measure of the amount of biological information transferred between groups, and McShane--Whitney widths are deterministic geometric compatibility ranges rather than probabilistic uncertainty intervals.

The sow-production application illustrates these mechanisms in a real longitudinal setting. Recent climatic and productive histories are represented as trajectories, farms define the group structure, and the complete procedure is evaluated under a causal month-by-month protocol. The empirical comparison shows that replacing exact local interpolation by robust latent fitting is accompanied by substantially narrower McShane--Whitney extension ranges and improved predictive performance, while the subsequent global--local construction provides a further improvement in the final evaluation period. These observations are intended to illustrate the behaviour of the mathematical framework rather than to establish general predictive superiority. A more extensive analysis of the production problem, including additional outcomes and a broader domain-specific interpretation, is reserved for a separate applied study.

The theoretical results remain conditional on a fixed pseudometric and a prescribed family of constrained pairs. The weighted trajectory pseudometric used in the application is therefore one particular realisation of the framework rather than an essential restriction of the theory. An important next step is to study the interaction between geometry learning and Lipschitz regression itself, including stability under perturbations of the learned metric and the joint estimation of trajectory geometry and latent response. Other natural directions include richer trajectory geometries, the behaviour of the shared regularity bound along its regularisation path, and the extension of the grouped construction to previously unseen groups. In this sense, the present work provides a link between the trajectory-space viewpoint developed in \citet{jonard2026} and a robust regression mechanism defined directly on that geometry, yielding an explicit route from trajectory representation to grouped and sequential prediction.

\section*{Acknowledgements and funding}

The authors would like to thank the Instituto de Ciencia y Tecnolog\'ia Animal (ICTA) at the Universidad Polit\'ecnica de Valencia for its assistance with the farm data, and Animal Data Analytics (ADA) for providing the dataset used in the empirical application presented in this study.

This research was funded by the Agencia Estatal de Investigaci\'on, grant PID2022-138342NB-I00, and by the European Union's Horizon Europe research and innovation programme under Grant Agreement No.~101059609 (Re-Livestock).

%\section*{Data and code availability}

\bibliographystyle{plainnat}
\bibliography{references}

\end{document}